\documentclass[number,citesort,seceqn,rotating,dvips]{arxbj}
\usepackage{upgreek}
\usepackage{graphicx}
\usepackage{accents}

\aid{0}
\volume{17}
\issue{3}
\pubyear{2011}
\firstpage{1063}
\lastpage{1094}
\doi{10.3150/10-BEJ298}

\makeatletter

\newtheorem{Lem}{Lemma}[section]
\newtheorem{Prop}{Proposition}[section]

\newcommand{\R}{\mathbb R}
\newcommand{\n}{^{(n)}}
\newcommand{\sirc}{\circ}

\newcommand{\tendL}{\stackrel{\mathcal{L}}{\longrightarrow}}
\newcommand{\varthetab}{\bolds{\vartheta}}
\newcommand{\Deltab}{\bolds{\Delta}}
\newcommand{\taub}{\bolds{\tau}}
\newcommand{\Gamb}{\bolds{\Gamma}}
\newcommand{\pr}{^{\prime}}
\newcommand{\utT}{\underaccent{\widetilde{}}{T}}
\newcommand{\utg}{\underaccent{\widetilde{}}{\gamma}}
\newcommand{\DS}[1]{#1}
\newcommand{\ny}{n\rightarrow\infty}

\makeatother

\begin{document}
\begin{frontmatter}

\title{A class of optimal tests for symmetry based on local Edgeworth
approximations}
\runtitle{Optimal tests for symmetry}

\begin{aug}
\author[a]{\fnms{Delphine} \snm{Cassart}\thanksref{a}},
\author[a,b,c]{\fnms{Marc} \snm{Hallin}\corref{}\thanksref{a,b,c}\ead[label=e1]{mhallin@ulb.ac.be}}
\and
\author[a,d]{\fnms{Davy} \snm{Paindaveine}\thanksref{a,d}}
\runauthor{D. Cassart, M. Hallin and D. Paindaveine}
\address[a]{Institut de Recherche en Statistique, ECARES, Universit\' e
Libre de Bruxelles, Brussels, Belgium}
\address[b]{Acad\' emie Royale de Belgique, Brussels, Belgium}
\address[c]{CentER, Tilburg University, Tilburg, The Netherlands}
\address[d]{D\' epartement de Math\' ematique, Universit\' e Libre de
Bruxelles, Brussels, Belgium}
\end{aug}

\received{\smonth{3} \syear{2007}}
\revised{\smonth{2} \syear{2010}}

%
\begin{abstract}
The objective of this paper is to provide, for the problem of
univariate symmetry (with respect to specified or unspecified
location), a concept of optimality, and to construct tests achieving
such optimality. This requires embedding symmetry into adequate
families of asymmetric (local) alternatives. We construct such
families by considering non-Gaussian generalizations of classical
first-order Edgeworth expansions indexed by a measure of skewness such
that (i) location, scale and skewness play well-separated roles
(diagonality of the corresponding information matrices) and (ii) the
classical tests based on the Pearson--Fisher coefficient of skewness
are optimal in the vicinity of Gaussian densities.
\end{abstract}

%

\end{frontmatter}

\section{Introduction}\label{introintro}

\subsection{Testing for symmetry}\label{introart1}
Symmetry is one of the most important and fundamental structural assumptions in
statistics, playing a major role, for instance, in the identifiability
of location or intercept under nonparametric conditions: see \cite{S1956,B1974,S1975}. This importance explains the
huge variety of existing testing procedures of the null hypothesis of
symmetry in an i.i.d.\ sample $X_1,\ldots, X_n$; see \cite{H1988}
for a survey.

Traditional tests of the null hypothesis of symmetry -- the hypothesis
under which $X_1-\theta\stackrel{d}{=} -(X_1-\theta)$ for some
location $\theta\in\R$, with $\stackrel{d}{=}$ standing for
equality in distribution -- are based on standardized empirical
third-order moments. Let $m\n_k(\theta):=n^{-1}\sum
_{i=1}^n(X_i-\theta)^k$ and $m\n_k :=m\n_k(\bar{X}\n)$, where
$\bar{X}\n:= n^{-1}\sum_{i=1}^nX_i$. When the location $\theta$ is
specified, the test statistic is
%
\begin{equation}\label{classicaltheta} S\n_1(\theta):=
n^{1/2}m\n_3(\theta)/ \bigl(m\n_6(\theta)\bigr)^{1/2},
\end{equation}
the null distribution of which, under finite sixth-order moments, is
asymptotically standard normal. When $\theta$ is unspecified, the
classical test is based on the empirical coefficient of skewness
%
\begin{equation}\label{b1} b_1\n:= m_3\n/s^3_n,
\end{equation}
where $s_n:=(m_2\n)^{1/2}$ stands for the empirical standard error in
a sample of size $n$. More precisely, this test relies on the
asymptotic standard normal distribution (still under finite moments of
order six) of
%
\begin{equation}\label{classical} S\n_2:=
n^{1/2}m_3\n\big/\bigl(m\n_6 - 6s_n^2m\n_4 + 9s_n^6\bigr)^{1/2},
\end{equation}
which, under Gaussian densities, asymptotically reduces to $\sqrt{n/6}
 b\n_1$.

These two tests are generally considered as Gaussian procedures,
although they do not require any Gaussian assumptions and despite the
fact that none of them can be considered optimal in any Gaussian sense,
since asymmetric alternatives clearly cannot belong to a Gaussian
universe. Despite the long history of the problem, the optimality
features of those classical procedures thus are all but clear, and
optimality issues, in that fundamental problem, remain essentially unexplored.

The main objective of this paper is to provide this classical testing
problem with a~concept of optimality that confirms practitioners'
intuition (i.e., justifying the $b\n_1$-based Gaussian practice),
and to construct tests achieving such optimality. This requires
embedding the null hypothesis of symmetry into adequate families of
asymmetric alternatives. We therefore define local (in the LeCam sense)
alternatives indexed by location, scale and a measure of skewness in
such a way that:
\begin{enumerate}[(ii)]
\item[(i)] Location, scale, and skewness play well-separated roles
(diagonality of the corresponding information matrices).
\item[(ii)] The traditional tests based on $b_1\n$ (more precisely,
based on $S\n_2$ given in (\ref{classical})) become locally and
asymptotically optimal in the vicinity of Gaussian densities.
\end{enumerate}

As we shall see, part (ii) of this objective is achieved 
by considering local first-order Edgeworth approximations of the form
%
\begin{equation}\label{EdgeG}
\phi(x-\theta) + n^{-1/2}\xi(x-\theta)\phi(x-\theta)\bigl((x-\theta
)^2-\kappa\bigr)
,
\end{equation}
where $\phi$ as usual stands for the standard normal density, $\kappa$
($=3$) is the Gaussian kurtosis coefficient, $\theta$ is a location
parameter, and $\xi$ a measure of skewness. Adequate modifications
of (\ref{EdgeG}), playing similar roles in the vicinity of
non-Gaussian standardized symmetric reference densities $f_1$, are
proposed in (\ref{modelart1}).

The resulting tests of symmetry (for specified as well as for
unspecified location $\theta$) are valid under a broad class of
symmetric densities, and parametrically efficient at the reference
(standardized) density $f_1$. Of particular interest are the
pseudo-Gaussian tests (associated with a Gaussian reference density;
see Proposition \ref{testGbart1}), which require finite moments of
order six and appear to be asymptotically equivalent (under their
\mbox{specified-$\theta$} version as well as under the unspecified-$\theta$
one) to the test (\ref{classical}) based on~$b\n_1$, and~the\vadjust{\eject}
Laplace tests (associated with a double-exponential reference density;
see Proposition~\ref{testGbart1}), which only require moments of order
four and are closely related with the tests against Fechner asymmetry
derived in \cite{CHP2008}.

These tests are of a parametric nature. Since the null hypothesis of
symmetry enjoys a~rich group invariance structure, classical maximal
invariance arguments naturally bring \textit{signs} and \textit{signed ranks}
into the picture. Such a nonparametric approach is adopted in a
companion paper \cite{CHP2010}, where we construct
signed-rank versions of the parametrically efficient tests proposed
here. These signed-rank tests are distribution-free (asymptotically so
in case of an unspecified location $\theta$) under the null hypothesis
of symmetry, and therefore remain valid under much milder
distributional assumptions (for the specified location case, they are
valid in the absence of \textit{any} distributional assumption).

The main technical tool throughout the paper is LeCam's asymptotic
theory of statistical experiments and the properties of locally
asymptotically normal (LAN) families. LAN has become a standard tool in
asymptotics: see \cite{lacamYang2000} or Chapters 6--9 of
\cite{vdV1988} for details. Log-likelihoods in a LAN family with
$k$-dimensional parameter $\varthetab$ admit local quadratic
approximations of the form $\taub\pr\Deltab\n_{\varthetab}- \frac
{1}{2}\taub\pr\Gamb_{\varthetab}\taub$, where the random vector
$\Deltab\n_{\varthetab}$, called a~\textit{central sequence}, is
asymptotically normal $\mathcal{N}(\Gamb_{\varthetab}\taub, \Gamb
_{\varthetab})$ under sequences of parameter values of the form
$\varthetab+ n^{-1/2}\taub$ (\textit{local alternatives}). Let $\phi
^*(\Deltab)$ be an optimal test (uniformly most powerful, maximin, most
stringent, \ldots) in the \textit{Gaussian shift model} describing a
hypothetical observation $\Deltab$ with distribution in the family $\{
\mathcal{N}(\Gamb_{\varthetab}\taub, \Gamb_{\varthetab}) \mid\taub \in\R^k
\}$ ($\Gamb_{\varthetab}$ specified). Those families are extremely
simple, and optimal tests in that context are well known (see, e.g., Section 11.9 of \cite{lecam1986}). Then, the sequence $\phi
^*(\Deltab\n_{\varthetab})$ is a sequence of \textit{locally
asymptotically optimal} (locally asymptotically \textit{uniformly most
powerful}, \textit{maximin}, \textit{most stringent}, \ldots) tests for
the original problem -- where optimality is based on the local
convergence of risk functions to the risk functions of Gaussian shift
experiments. More analytical characterizations can be found, for
instance, in \cite{CHS1996}.

\subsection{Outline of the paper}

The problem we are considering throughout is that of testing the null
hypothesis of symmetry. In the notation of Section \ref{introart1},
$\xi$ (see (\ref{modelart1}) for a more precise definition) is thus
the parameter of interest; the location $\theta$ and the standardized
null symmetric density~$f_1$ either are specified or play the role of
nuisance parameters, whereas the scale~$\sigma$ (not necessarily a
standard error) always is a nuisance.

The paper is organized as follows. In Section \ref{introEdg} we
describe the Edgeworth-type families of local alternatives, extending
(\ref{EdgeG}), that we are considering. Section \ref{ULAN}
establishes the local and asymptotic normality (with respect to
location, scale and the asymmetry parameters) result to be used
throughout; actually, we establish a slightly stronger version
of LAN, called ULAN (uniform LAN), which allows us to handle the
problems related with estimated nuisance parameters. The classical
LeCam theory then is used in Section \ref{Sec3.2} for developing
asymptotically optimal procedures for testing symmetry ($\xi=0$), with
specified or unspecified location $\theta$ but specified standardized
symmetric density $f_1$.
The more realistic case of an unspecified $f_1$ is treated in Section
\ref{parametriart1}, where we obtain versions of the optimal (at given
$f_1$) tests that remain valid under $g_1 \ne f_1$, for specified
(Section \ref{321}) and unspecified (Section \ref{322}) location
$\theta$, respectively. The particular case of pseudo-Gaussian
procedures (which are optimal for Gaussian $f_1$ but valid under any
symmetric density with finite moments of order six) is studied in
detail in Section \ref{testspseudogauss} and their relation with
classical tests of symmetry is discussed. We also show that the Laplace
tests (which are optimal for double-exponential $f_1$ but valid under
any symmetric density with finite fourth-order moment) are closely
related to the Fechner-type tests derived in \cite{CHP2008}.
The finite-sample performances of these tests are investigated via
simulations in Section \ref{simu}, where they are applied to the
classical skew-normal and skew-$t$ densities.

%

\section{A class of locally asymptotically normal families of
asymmetric distributions }\label{Sec2}

\subsection{Families of asymmetric densities based on Edgeworth
approximations}\label{introEdg}
Denote by $\mathbf{X}^{(n)}:=(X_1^{(n)}, \ldots, X_n^{(n)}),\ n \in
\mathbb N$ an $\mbox{i.i.d.}$ $n$-tuple of observations with common
density~$f$. The null hypotheses we are interested in are:
\begin{enumerate}[(a)]
\item[(a)]The hypothesis $\mathcal{H}\n_\theta$ of symmetry with
respect to specified location $\theta\in\R$: un\-der~$\mathcal{H}\n
_\theta$, the $X_i$'s have density function\vspace*{-3pt}
%
\begin{equation}\label{fsym}
x\mapsto f(x):= {\sigma}^{-1}f_1\bigl(({x-\theta})/{\sigma}\bigr)\vspace*{-3pt}
\end{equation}
(all densities are over the real line, with respect to the Lebesgue
measure) for some unspecified $\sigma\in\R^+_0$, where $f_1$ belongs
to the class of standardized symmetric densities\vspace*{-3pt}
\[
\mathcal{F}_0 := \biggl\{h_1\dvtx h_1(-z)=h_1(z) \mbox{ for all $z\in\R$ and } \int_{-\infty}^1 h_1(z) \,\mathrm{d}z
= 0.75 \biggr\}.\vspace*{-3pt}
\]
The scale parameter $\sigma$ (associated with the symmetric density
$f$) that we are considering here thus is not the standard
deviation, but the median of the absolute deviations $\vert
X_i-\theta\vert$; this avoids making any moment assumptions.
\item[(b)]The hypothesis $\mathcal{H}\n:=\bigcup_{\theta\in\R
}\mathcal{H}\n_\theta$ of symmetry with respect to unspecified
location (and scale): there exist $(\theta, \sigma)$ such that the
$X_i$'s have density (\ref{fsym}).
\end{enumerate}
In both cases, the standardized density $f_1$ may be specified (Section
\ref{Sec3.2}) or not (Sections~\ref{parametriart1}--\ref{testsLaplace}). The specified-$f_1$ problem,
however, mainly serves as a preparation for the more realistic
unspecified-$f_1$ one.

As explained in the introduction, a characterization of efficient
testing requires the definition of families of asymmetric alternatives
exhibiting some adequate structure, such as local asymptotic normality,
at the null. For a selected class of densities $f$ enjoying the
required regularity assumptions, we therefore are embedding the null
hypothesis of symmetry into families of distributions indexed by
$\theta\in\R$ (location), $\sigma\in\R^+_0$ (scale) and a
parameter $\xi\in\R$ characterizing asymmetry. More precisely,
consider the class $\mathcal{F}_1$ of densities $f_1$ satisfying:
\begin{enumerate}[(iii)]
\item[(i)] (symmetry and standardization) $f_1\in \mathcal{F}_0$;
\item[(ii)] (absolute continuity) there exists $\dot{f}_1$ such that,
for all $z_1<z_2$,\vspace*{-3pt}
\[
f_1(z_2) - f_1(z_1)=\int_{z_1}^{z_2}\dot
{f}_1(z) \,\mathrm{d}z;\vspace*{-3pt}
\]
\item[(iii)] (strong unimodality) $\DS{z\mapsto \phi_{f_1}(z):=
{-\dot{f}_1(z)}/{f_1(z)}}$ is monotone increasing;
\item[(iv)] (finite Fisher information) $\DS{ \mathcal{K} (f_1) :=
\int_{- \infty}^{+ \infty} z^4 \phi^2_{f_1}(z) f_1(z) \,\mathrm{d}z}$, hence
also, under strong unimodality,\vspace*{-3pt}
\[
\mathcal{I} (f_1) := \int_{-
\infty}^{+ \infty} \phi^2_{f_1}(z) f_1(z) \,\mathrm{d}z \quad  \mbox{and}
\quad \mathcal{J} (f_1) := \int_{- \infty}^{+ \infty} z^2 \phi
^2_{f_1}(z) f_1(z)\,\mathrm{d}z\vspace*{-3pt}
\]
are finite;
\item[(v)] (polynomial tails)\vspace*{-3pt}
\[
\int_y^{\infty}
f_1(z) \,\mathrm{d}z = \mathrm{O}(y^{-\beta})\qquad  \mbox{as }y\to\infty\vspace*{-3pt}
\]
for some $\beta>0$ and\vspace*{-3pt}
\[
\phi_{f_1}(z)=
\mathrm{o}(z^{\beta/2-2}) \qquad  \mbox{as }z \rightarrow\infty.\vspace*{-3pt}
\]
\end{enumerate}
That class
$\mathcal{F}_1$ thus consists of all symmetric standardized densities
$f_1$ that are absolutely continuous, strongly unimodal (that is,
log-concave) and have finite information $\mathcal{I }(f_1)$ and
$\mathcal{J }(f_1)$ for location and scale, and, as we shall see,
$\mathcal{K} (f_1)$ for asymmetry, with tails satisfying (v).

For all $f_1\in\mathcal{F}_1$, denote by $\kappa(f_1):=\mathcal
{J}(f_1)/\mathcal{I}(f_1)$ the ratio of information for scale and
information for location; $\kappa(f_1)$
, as we shall see, for Gaussian density ($f_1=\phi_1$) reduces to
kurtosis ($\kappa(\phi_1)=3$), and can be interpreted as a \textit{generalized kurtosis coefficient}. Finally, write ${\mathrm{P}}\n_{\theta,
\sigma, \xi; f_1}$ for the probability distribution of $
\mathbf{X}^{(n)}$ when the $X_i$'s are \mbox{i.i.d.} with density\vspace*{-3pt}
%
\begin{eqnarray}\label{modelart1}
f(x)&=&\sigma^{-1} f_1 \biggl(\frac{x-\theta}{\sigma} \biggr) - \xi\sigma
^{-1} {\dot{f}}_1 \biggl(\frac{x-\theta}{\sigma} \biggr) \biggl( \biggl(\frac{x-\theta
}{\sigma} \biggr)^2-\kappa(f_1) \biggr)I[\vert x -\theta\vert\leq\sigma\vert
z^* \vert]
\nonumber\\[-10pt]\\[-10pt]
&& - \operatorname{sign}(\xi) \sigma^{-1} f_1 \biggl(\frac{x-\theta}{\sigma} \biggr) \{
I[ x-\theta> \operatorname{sign}(-\xi) \sigma\vert z^*\vert ]- I[ x-\theta<
\operatorname{sign}(\xi)\sigma\vert z^*\vert ] \}.\nonumber\hspace*{25pt}\vspace*{-3pt}
\end{eqnarray}
Here $\theta\in\R$ and $\sigma\in\R^+$ clearly are location and
scale parameters, $\xi\in\R$ is a measure of skewness, $\kappa
(f_1)$ (strictly positive for $f_1\in\mathcal{F}_1$) the generalized
kurtosis coefficient just defined and $z^*$ the unique (for $\xi$
small enough; unicity follows from the monotonicity of~$\phi_{f_1} $)
solution of $f_1(z^*) = \xi\dot{f}_1(z^*) ((z^*)^2 - \kappa(f_1) )$.
The function $f$ defined in (\ref{modelart1}) is indeed a probability
density (non-negative, integrating up to one), since it is obtained by
adding and subtracting the same probability mass
\[
\frac{\vert\xi\vert}{\sigma}\int_\theta^\infty\min
\biggl({\dot{f}}_1 \biggl(\frac{x-\theta}{\sigma} \biggr) \biggl( \biggl(\frac{x-\theta}{\sigma
} \biggr)^2-\kappa(f_1) \biggr) , f_1 \biggl(\frac{x-\theta}{\sigma} \biggr) \biggr)\,\mathrm{d}x
\]
on both sides of $\theta$ (according to the sign of $\xi$). Note that
$\xi>0$ implies $f(x)=0 $ for $x-\theta< -\sigma\vert z^*\vert$ and
$f(x)= {2}{\sigma}^{-1}f_1((x-\theta)/\sigma) $ for $x-\theta>
\sigma\vert z^*\vert$. Moreover, $x\mapsto f(x)$ is continuous
whenever $\dot{f}_1(x)$ is; vanishes for $x \leq\theta+ \sigma z^*$
if $\xi>0$, for $x \geq\theta+ \sigma z^*$ if $\xi<0$; and is left-
or right-skewed according as $\xi<0$ or $\xi>0$. As for $z^*$, it
tends to $-\infty$ as $\xi\downarrow0$, to $\infty$ as $\xi
\uparrow0$; in the Gaussian case, it is easy to check that $\vert
z^*\vert=\mathrm{O}(\vert\xi\vert ^{-1/3})$ as $\xi\rightarrow0$.

The intuition behind this class of alternatives is that, in the
Gaussian case, (\ref{modelart1}), with $\xi= n^{-1/2}\tau$ yields
(for $x\in[\theta\pm\sigma z^*]$) the first-order Edgeworth
development of the density of the standardized mean of an \mbox
{i.i.d.} $n$-tuple of variables with third-order moment $6\tau\sigma
^3$ (where standardization is based on the median $\sigma$ of absolute
deviations from $\theta$). For a~``small'' value of the asymmetry
parameter $\xi$, of the form $ n^{-1/2}\tau$, (\ref{modelart1}) thus
describes the type of deviation from symmetry that corresponds to the
classical central limit context. Hence, if a Gaussian density is
justified as resulting from the additive combination of a~large number
of small independent symmetric shocks, the locally asymmetric $f$
results from the same additive combination of independent but slightly
skew shocks. As we shall see, the locally optimal test in such a case
happens to be the traditional test based on~$b_1\n$ (see (\ref{b1})).

Besides the Gaussian one (with standardized density $\phi_1(z):=\sqrt {
{a}/{2\uppi}}\exp(-az^2/2)$), interesting special cases of
(\ref{modelart1}) are obtained in the vicinity of:
\begin{enumerate}[(iii)]
\item[(i)] The double-exponential or Laplace distributions, with standardized density\vspace*{-3pt}
\[
f_1(z)=f_\mathcal{L}(z):= ({1}/{2d}) \exp(- {|z|}/{d}),\vspace*{-3pt}
\]
$\mathcal{I}(f_1)={1}/{d^2}$, $\mathcal{J}(f_1)=2$ and $\mathcal
{K}(f_1)=24d^2$.
\item[(ii)] The logistic distributions, with standardized density\vspace*{-3pt}
\[
f_1(z)=f_{\mathrm{Log}}(z):=\sqrt{b} \exp\bigl(-\sqrt{b}z\bigr)/\bigl(1+\exp\bigl(-\sqrt{b}z\bigr)\bigr)^2,
\vspace*{-3pt}
\]
$\mathcal{I}(f_1)= {b}/{3}$, $\mathcal{J}(f_1)= {(12+\uppi^2)}/{9}$
and $\mathcal{K}(f_1)={\uppi^2} (120+7\uppi^2)/ 45b$.
\item[(iii)] The power-exponential distributions, with standardized densities\vspace*{-3pt}
\[
f_1(z)=f_{\exp_{\eta}}(z):= C_{\exp_{\eta}}\exp
(-(g_\eta z)^{2\eta}),\vspace*{-3pt}
\]
$\eta\in\mathbb N_0$, $\mathcal{I}(f_1)= 2g_\eta^2\eta\Gamma
(2-1/2\eta)/\Gamma(1+1/2\eta)$, $\mathcal{J}(f_1)=1+2\eta$ and
$\mathcal{K}(f_1)= {2g_\eta\eta}/\break{\Gamma(1+1/2\eta)}$
(the positive constants $C_{\exp_{\eta}}$, $a$, $b$, $d$ and
$g_\eta$ are such that $f_1\in\mathcal{F}_1$).
\end{enumerate}

Although not strongly unimodal, the Student distributions with $\nu>2$
degrees of freedom also can be considered here (strong unimodality
indeed is essentially used as a~sufficient condition for the existence
of $z^*$ in (\ref{modelart1}) -- an existence that can be checked
directly in the Student case). Standardized Student densities take the form\vspace*{-3pt}
\[
f_1(z)=f_{t_{\nu}}(z):=C_{t_{\nu}}(1+a_{\nu}z^2/\nu)^{-(\nu+1)/2},\vspace*{-3pt}
\]
with $\mathcal{I}(f_1)=a_{\nu} {(\nu+1)}/{(\nu+3)}$, $ \mathcal
{J}(f_1)=3 {(\nu+1)}/{(\nu+3)}$ and $\mathcal{K}(f_1)= {15} \nu
(\nu+1)/a_{\nu}(\nu-2)(\nu+3)$ ($C_{t_{\nu}}$ and $a_{\nu}$
are normalizing constants). Note that the corresponding Gaussian
values, namely $\mathcal{I}(\phi_1)=a \approx 0.4549$, $\mathcal
{J}(\phi_1)=3$ and $\mathcal{K}(\phi_1)= {15}/{a}$, are obtained by
taking limits as $\nu\rightarrow\infty$.

\begin{figure}

\includegraphics{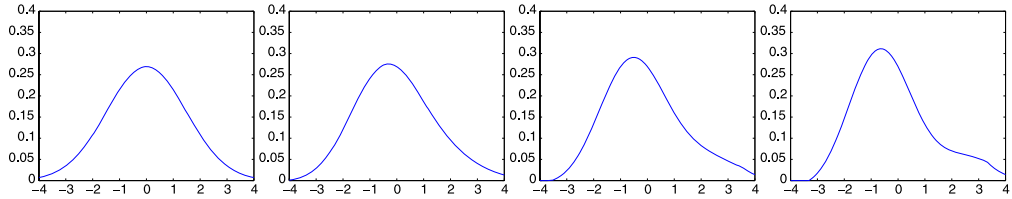}
\vspace*{-6pt}
\caption{Graphical representation of the Gaussian Edgeworth
family (\protect\ref{modelart1}) ($f_1=\phi_1$), for $\xi=0$, $0.05$, $0.10$
and $0.15$.}\label{graphmodart1}
\vspace*{-5pt}
\end{figure}

\begin{figure}[b]

\includegraphics{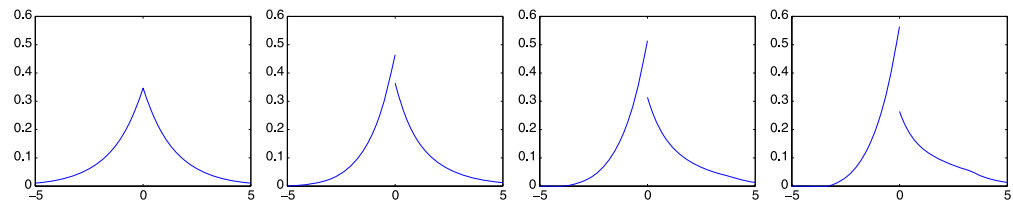}%
\vspace*{-6pt}
\caption{Graphical representation of the double-exponential
Edgeworth family (\protect\ref{modelart1}) ($f_1=f_\mathcal{L}$), for $\xi
=0$, $0.05$, $0.10$ and $0.15$.}
\label{graph2modart1}
\vspace*{-5pt}
\end{figure}

Figures \ref{graphmodart1} and \ref{graph2modart1} provide graphical
representations of some densities in the Gaussian ($f_1=\phi_1$) and
double-exponential ($ f_1=f_\mathcal{L}$) Edgeworth families (\ref
{modelart1}), respectively. In the Gaussian case, the skewed densities
are continuous, while the double-exponential ones, due to the
discontinuity of $\dot{f}_\mathcal{L}(x)$ at $x=0$, exhibit a
discontinuity at the origin.

\subsection{Uniform local asymptotic normality (ULAN)}\label{ULAN}
The main technical tool in our derivation of optimal tests is the
uniform local asymptotic normality (ULAN), with respect to $\bolds{\vartheta}:=(\theta, \sigma, \xi)\pr$, at $(\theta, \sigma
,0)\pr$, of the parametric families\vspace*{-3pt}
%
\begin{equation}\label{famparam}
\mathcal{P}_{f_1}^{(n)}:= \bigcup_{\sigma> 0}\mathcal{P}_{\sigma;f_1}^{(n)}
:= \bigcup_{\sigma> 0} \bigl\{{\mathrm{P}}_{\theta, \sigma, \xi;f_1}^{(n)}
\mid
\theta\in\mathbb R,  \xi\in\mathbb R \bigr\},\vspace*{-5pt}
\end{equation}
where $f_1\in\mathcal{F}_1$. More precisely, the following result
holds (see the \hyperref[appendixart1]{Appendix} for proof).

\begin{Prop}[(ULAN)] \label{lanart1}
For any $f_1\in\mathcal{F}_1$, $\theta\in\R$, and $\sigma\in\R
^+_0$, the family $\mathcal{P}_{f_1}^{(n)}$ is ULAN at $(\theta,
\sigma,0)\pr$, with (writing $Z_i$ for $Z_i^{(n)}(\theta, \sigma
):=\sigma^{-1}(X_i^{(n)} - \theta)$ and $\phi_{f_1}$ for $-\dot
{f_1}/f_1$) central sequence\vspace*{-3pt}
%
\begin{eqnarray} \label{centralseqart1}
\bolds{\Delta}_{f_1}^{(n)}(\bolds{\vartheta})&=:& \pmatrix{
{\Delta}_{f_1;1}^{(n)}(\bolds{\vartheta}) \cr
{\Delta}_{f_1;2}^{(n)}(\bolds{\vartheta})\cr
{\Delta}_{f_1;3}^{(n)}(\bolds{\vartheta})
}\nonumber\\[-8pt]\\[-8pt]
&=\hspace*{3pt}& n^{-1/2} \sum_{i=1}^n \pmatrix{
\sigma^{-1}\phi_{f_1}(Z_i) \cr
\sigma^{-1}\bigl(\phi_{f_1}(Z_i) Z_i - 1\bigr)\cr
\phi_{f_1}(Z_i) \bigl( Z_i^2 - \kappa(f_1) \bigr)
}\nonumber\vspace*{-3pt}
\end{eqnarray}
and full-rank information matrix
%
\begin{equation}\label{infomatart1}
\bolds{\Gamma}_{f_1}(\bolds{\vartheta}) = \pmatrix{
\sigma^{-2}\mathcal{I}(f_1) & 0 & 0\cr
0 & \sigma^{-2}\bigl(\mathcal{J}(f_1)-1\bigr) & 0 \cr
0 & 0 & \gamma(f_1)
},
\end{equation}
where $\gamma(f_1):=\mathcal{K}(f_1)- {\mathcal{J}^2(f_1)}/{\mathcal
{I}(f_1)}$.

More precisely, for any $\bolds{\vartheta}^{(n)}:=(\theta^{(n)},
\sigma^{(n)}, 0)\pr$ such that $\theta\n-\theta=\mathrm{O}(n^{-1/2})$ and
$\sigma\n- \sigma=\mathrm{O}(n^{-1/2})$, and for any bounded sequence $\bolds{\tau}^{(n)}=(t^{(n)}, s^{(n)}, \tau^{(n)})\pr \in\mathbb R^3$, we
have, under ${\mathrm{P}}_{\bolds{\vartheta}^{(n)};f_1}^{(n)}$, as $n
\rightarrow\infty$,
\begin{eqnarray*}
\Lambda_{\bolds{\vartheta}^{(n)}+n^{-1/2}\bolds{\tau}^{(n)} / \bolds{\vartheta}^{(n)};f_1}^{(n)} &:=&\log\biggl( \frac{\mathrm{d}{\mathrm{P}}_{\bolds{\vartheta}^{(n)}+n^{-1/2}\bolds{\tau}^{(n)};f_1}^{(n)}}{\mathrm{d}{\mathrm{P}}_{\bolds{\vartheta}^{(n)};f_1}^{(n)}} \biggr)\\
&=&\bolds{\tau}^{(n)\prime} \bolds{\Delta}_{f_1}^{(n)}\bigl(\bolds{\vartheta}\n\bigr) - \frac{1}{2} \bolds{\tau}^{(n)\prime} \bolds{\Gamma
}_{f_1}(\bolds{\vartheta}) \bolds{\tau}^{(n)} +\mathrm{o}_{\mathrm{P}}(1)
\end{eqnarray*}
and
\[
{\bolds{\Delta}_{f_1}^{(n)}\bigl( \bolds{\vartheta}^{(n)} \bigr) \stackrel
{\mathcal{L}}{\longrightarrow} \mathcal{N}(\bolds{0}, \bolds{\Gamma
}_{f_1}(\bolds{\vartheta}))}.
\]
\end{Prop}

The diagonal form of the information matrix $\bolds{\Gamma}_{f_1}(\bolds{\vartheta})$ confirms that location, scale and skewness, in the
parametric family (\ref{famparam}), play distinct and well-separated
roles. The practical consequences of that orthogonality (in the sense
of information) property, as we shall see, are twofold:
\begin{enumerate}[(b)]
\item[(a)] The fact that location and scale are unspecified has no
cost in terms of efficiency and under specified standardized density $f_1$ when
testing for symmetry.
\item[(b)] Substituting root-$n$-consistent (and, in principle, duly
\textit{discretized}: see assumption (C2) below) estimators for the true
values has no impact on the asymptotic validity and local powers, under
specified
standardized density $f_1$, of tests for symmetry.
\end{enumerate}
Note that orthogonality between the scale and skewness components of
$\bolds{\Delta}_{f _1}^{(n)}(\bolds{\vartheta})$ automatically follows
from the symmetry of $f_1$, while for location and skewness, this
orthogonality is a consequence of the definition of $\kappa(f_1)$. The
Gaussian versions of (\ref{centralseqart1}) and (\ref{infomatart1}) are
\[
\bolds{\Delta}_{\phi_1}^{(n)}(\bolds{\vartheta})= n^{-1/2}\sum
_{i=1}^n \pmatrix{
a\sigma^{-1}Z_i\cr
\sigma^{-1}(aZ_i^2 - 1)\cr
\displaystyle aZ_i\biggl(Z_i^2-\frac{3}{a}\biggr)
}
\]
and
\[
\bolds{\Gamma}_{\phi_1}(\bolds{\vartheta}) = \pmatrix{
a\sigma^{-2} & 0 & 0\cr
0 &2 \sigma^{-2} & 0 \cr
0 & 0 &6/a
},
\]
respectively (recall that $a \approx 0.4549$).

\section{Locally asymptotically optimal tests}\label{Sec3}

The various test statistics described in this section are listed, for
easy reference, in Table~\ref{summary}, placed at the end of Section \ref{conclu}.

\subsection{Locally asymptotically optimal tests: Specified
density}\label{Sec3.2}
For specified $f_1\in\mathcal{F}_1$, consider the null hypothesis
\[
\mathcal{H}\n_{\theta;f_1}:=\bigcup_{\sigma\in\R^+_0}\bigl\{{\mathrm{P}}\n
_{\theta,\sigma, 0;f_1}\bigr\}
\]
of symmetry with respect to some specified location $\theta$, and the
null hypothesis
\[
\mathcal{H}\n_{f_1}:=\bigcup_{\theta\in\R}\bigcup_{\sigma\in\R
^+_0}\bigl\{{\mathrm{P}}\n_{\theta,\sigma, 0;f_1}\bigr\}
\]
of symmetry with respect to unspecified $\theta$. ULAN and the
diagonal structure of (\ref{infomatart1}) imply that substituting
discretized root-$n$-consistent estimators $\hat{\theta}$ and $\hat
{\sigma}$ for the unknown~$\theta$ and $\sigma$ has no influence,
asymptotically, on the $\xi$-part of the central sequence.

Recall that a sequence of estimators $\hat\lambda\n$ defined in a
sequence of experiments $\{{\mathrm{P}}\n_\lambda\mid\lambda\in\Lambda
\}$ indexed by some parameter $\lambda$ is \textit{root-$n$-consistent}
and \textit{asymptotically discrete} if, under ${\mathrm{P}}\n_\lambda$, as
$\ny$:
\begin{enumerate}[(C2)]
\item[(C1)] $\hat\lambda\n- \lambda= \mathrm{O}_{\mathrm{P}}(n^{-1/2})$.

\item[(C2)] The number of possible values of $\hat\lambda\n$ in
balls with $\mathrm{O}(n^{-1/2})$ radius centered at~$\lambda$ is bounded as
$\ny$.
\end{enumerate}
An estimator $\lambda\n$ satisfying (C1) but not (C2) is easily
discretized by letting, for some arbitrary constant $c>0$, $\lambda\n
_{\#}:=(cn^{1/2})^{-1}\operatorname{sign}(\lambda\n)\lceil cn^{1/2}\vert
\lambda\n\vert\rceil$, which satisfies both (C1) and (C2).
Subscripts $_\#$ in the sequel are used for estimators ($\hat\theta_\#
$, $\hat\sigma_\#$, \ldots) satisfying (C1) and (C2).
It should be noted, however, that (C2) has no implications in practice,
where $n$ is fixed, as the discretization constant $c$ can be chosen
arbitrarily large.

It follows from the diagonal form of the information matrix (\ref
{infomatart1}) that locally uniformly asymptotically most powerful
tests of $\mathcal{H}\n_{\theta;f_1}$ (resp., of $\mathcal{H}\n
_{f_1}$) can be based on ${\Delta}_{f_1;3}^{(n)}( {\theta}, \hat
{\sigma}_\#,0)$ (resp., on ${\Delta}_{f_1;3}^{(n)}( \hat{\theta}_\#
, \hat{\sigma}_\#,0)$), hence on $T^{(n)}_{f_1}( {\theta}, \hat
{\sigma}_\# )$ (resp., on $T^{(n)}_{f_1}( \hat{\theta}_\#, \hat
{\sigma}_\#)$), where
%
\begin{equation}\label{stat1art1}
T^{(n)}_{f_1} ( {\theta}, {\sigma}):= \frac{1}{\sqrt{n\gamma
(f_1)}}\sum_{i=1}^{n} \phi_{f_1}( Z_i( {\theta}, {\sigma}) ) \bigl(
Z_i^2( {\theta} , {\sigma} ) - \kappa(f_1) \bigr).
\end{equation}
Root-$n$-consistent (under the null hypothesis of symmetry) estimators
of $\theta$ and $\sigma$ that do not require any moment assumptions
are, for instance, the medians $\hat\theta:=\operatorname{Med}(X_i^{(n)})$ and
$\hat\sigma:=\operatorname{Med}(\vert X_i^{(n)} - \hat\theta\vert)$ of the
$X_i^{(n)}$'s and of their absolute deviations from $\hat\theta$,
respectively.

The following proposition then results from classical results on ULAN
families (see, e.g., Chapter 11 of \cite{lecam1986}).

\begin{Prop}\label{test1art1}
Let $f_1\in\mathcal{F}_1$. Then:
\begin{enumerate}[(ii)]
\item[(i)] $T^{(n)}_{f_1}( \hat{\theta}_\#, \hat{\sigma}_\#
)=T^{(n)}_{f_1}( {\theta}, {\sigma} )+\mathrm{o}_{\mathrm{P}}(1)$ is asymptotically
normal, with mean zero under ${\mathrm{P}}\n_{\theta, \sigma, 0;f_1}
$, mean $\tau \gamma^{1/2}(f_1)$ under $
{\mathrm{P}}_{\theta,\sigma,n^{-1/2}\tau; f_1}^{(n)}
$ and variance one under both.
\item[(ii)] The sequence of tests rejecting the null hypothesis of
symmetry (with standardized density $f_1$) whenever
$T^{(n)}_{f_1}( {\theta}, \hat{\sigma}_\#)$ (resp., $T^{(n)}_{f_1}(
\hat{\theta}_\#, \hat{\sigma}_\#)$) exceeds the $(1-\alpha)$
standard normal quantile $z_\alpha$ is locally asymptotically most
powerful at asymptotic level $\alpha$ for $\mathcal{H}\n_{\theta
;f_1}$ (resp., for $\mathcal{H}\n_{f_1}$) against $\bigcup_{\xi>0}
\bigcup_{ \sigma\in\R^+_0}\{{\mathrm{P}}\n_{\theta, \sigma, \xi
;f_1}\}$ (resp., $\bigcup_{\xi>0}\bigcup_{\theta\in\R}\bigcup_{
\sigma\in\R^+_0}\{{\mathrm{P}}\n_{\theta, \sigma, \xi;f_1}\}$).
\end{enumerate}
\end{Prop}

This confirms that unspecified location $\theta$ and scale $\sigma$
do not induce any loss of efficiency when the standardized density
$f_1$ itself is specified.

The Gaussian version of (\ref{stat1art1}) is
\[
T^{(n)}_{\phi_1} ( {\theta}, {\sigma}):= \sqrt{\frac{a^3}{6n}}\sum
_{i=1}^{n} Z_i( {\theta}, {\sigma}) \biggl( Z_i^2( {\theta} , {\sigma} )
- \frac{3}{a} \biggr) = \sqrt{\frac{a}{6n}}\sum_{i=1}^{n} \bigl( aZ_i^3(
{\theta}, {\sigma}) - 3Z_i( {\theta} , {\sigma} ) \bigr) ;
\]
thanks to the linearity of Gaussian scores, it easily follows from a
traditional Slutsky argument that $\hat\theta$ and $\hat\sigma$ in
$T^{(n)}_{\phi_1} ( \hat\theta, \hat\sigma)$ need not be discretized.
Under Gaussian densities, both $T^{(n)}_{\phi_1} ( {\hat\theta},
{\hat\sigma})$ and $T^{(n)}_{\phi_1} ( {\theta}, {\hat\sigma})$
are asymptotically equivalent to $T^{(n)}_{\phi_1} ( \bar{X}^{(n)},
{\hat\sigma})=(na^3/6)^{1/2} m\n_3/\hat\sigma^3=\sqrt{n/6} \ b\n
_1 +\mathrm{o}_{\mathrm{P}}(1)$, that is, to $S\n_2$ given in (\ref{classical}).
The latter is thus locally asymptotically optimal under Gaussian
assumptions, whether $\theta$ is specified
or not, whereas the specified-$\theta$ test based on $m\n_3(\theta
)/(m\n_6(\theta))^{1/2}$ (more precisely, on $S\n_1(\theta)$ given
in (\ref{classicaltheta})) is suboptimal. The fact that $m\n_3(\hat
\theta)$ yields a better performance than $m\n_3(\theta)$ under
specified location $\theta$ (see the comments after Proposition \ref
{testGbart2} for a~comparison of local powers) looks puzzling at first
sight. The reason is that orthogonality, in the Fisher information
sense, between asymmetry and location, is a ``built-in'' feature of
Edgeworth families. Contrary to $S\n_2$, which is shift invariant,
$m\n_3(\theta)$ and $S\n_1(\theta)$ are sensitive to location
shifts. Therefore, tests based on $S\n_1(\theta)$ are ``wasting away''
some power on location alternatives (which are irrelevant when $\theta
$ is specified), to the detriment of asymmetry alternatives.

Locally asymptotically maximin two-sided tests are easily derived along
the same lines.

\subsection{Locally asymptotically optimal tests: Unspecified
density}\label{parametriart1}
The tests based on (\ref{stat1art1}) achieve local and asymptotic
optimality at correctly spe\-ci\-fied~$f_1$, which sets the parametric
efficiency bounds for the problem, but has limited practical value, as
these tests are not valid anymore under density $g_1 \neq f_1$. If
Proposition \ref{test1art1} is to be adapted to the more realistic
null hypotheses $\mathcal{H}\n_{\theta}:=\bigcup_{g_1}\mathcal
{H}\n_{\theta;g_1}$ and $\mathcal{H}\n:=\bigcup_{g_1}\mathcal
{H}\n_{ g_1}$ under which the (symmetric) density remains unspecified,
the test statistic~$T^{(n)}_{f_1}$ needs to be adapted in order to cope
with three problems:
its centering under the null and $g_1\neq f_1$ ($T^{(n)}_{f_1}$ in the
previous section only had to be centered under density $g_1=f_1$), its
scaling (same remark) and (still under the null and $g_1\neq f_1$)
the impact of the substitution of estimators $\hat\sigma$ (and $\hat
\theta$) for the unspecified values of $\sigma$ (and $\theta$) on
its asymptotic distribution.

\subsubsection{Specified location}\label{321}

Let us first assume that both $\theta$
and $\sigma$ are specified. The test statistic $T^{(n)}_{f_1}$ in
Proposition \ref{test1art1} is essentially a scaled version of $\Delta
_{f_1;3}^{(n)}$. More generally, write $\Delta_{f_1;3}^{(n)}(\kappa)$
for $ n^{-1/2}\sum_{i=1}^n \phi_{f_1}(Z_i) ( Z_i^2 - \kappa) $,
where $\kappa\in\R^+_0$ denotes a strictly positive real number; an
adequate sample-based value will be selected later on. Note
that $\Delta_{f_1;3}^{(n)}(\kappa)$ remains centered under ${\mathrm{P}}_{\theta, \sigma, 0;g_1}^{(n)}$, irrespective of the choice of
$\kappa$. Indeed, the functions $z\mapsto\phi_{f_1}(z)z^2$ and
$z\mapsto\phi_{f_1}(z) $ are skew symmetric, and their expectations
under any symmetric density are automatically zero -- provided that
they exist. The variance under ${\mathrm{P}}_{\theta, \sigma,
0;g_1}^{(n)}$of $\Delta_{f_1;3}^{(n)}(\kappa)$ is then\looseness=-1\vspace*{-2pt}
\[
\gamma^\kappa_{g_1}(f_1):=
\mathrm{E}_{g_1}\bigl[ \bigl( \phi_{f_1}(Z_i) (Z_i^2-\kappa)
\bigr)^2\bigr]
=
\mathcal{K}_{g_1}(f_1) - 2\kappa\mathcal{J}_{g_1}(f_1 ) + \kappa^2
\mathcal{I}_{g_1}(f_1),\vspace*{-2pt}
\]
where
\[
\mathcal{I}_{g_1}(f_1):= \int_{-\infty}^{\infty}\phi^2_{f_1}(z)
g_1(z) \,\mathrm{d}z,\qquad
\mathcal{J}_{g_1}(f_1 ):=\int_{-\infty}^{\infty}z^2\phi
^2_{f_1}(z)g_1(z) \,\mathrm{d}z
\]
and (still, provided that those integrals exist)\vspace*{-2pt}
\[
\mathcal{K}_{g_1}(f_1):= \int_{- \infty}^{+ \infty} z^4 \phi
_{f_1}^2 (z) g_1(z) \,\mathrm{d}z .\vspace*{-2pt}
\]

We know from LeCam's third lemma that, under ${\mathrm{P}}_{\theta, \sigma
, 0;g_1}^{(n)}$, the impact on $\Delta_{f_1;3}^{(n)}(\kappa)$ of an
estimated scale depends on the asymptotic joint distribution (under
${\mathrm{P}}_{\theta, \sigma, 0;g_1}^{(n)}$) of $\Delta
_{f_1;3}^{(n)}(\kappa)$ and $\Delta\n_{g_1;2}$. More precisely,
LeCam's third lemma (see, e.g., page 90 of \cite{vdV1988}) tells
us that if, under ${\mathrm{P}}_{\theta, \sigma, 0;g_1}^{(n)}$,
%
\begin{equation}  \label{d11}
\pmatrix{
\Delta_{f_1;3}^{(n)}(\kappa)\cr
\displaystyle \log\biggl(\frac{\mathrm{d}{\mathrm{P}}\n_{\theta, \sigma+ n^{-1/2}\tau, 0 ;g_1}}
{\mathrm{d}{\mathrm{P}}\n_{\theta, \sigma, 0 ;g_1}}
\biggr)
}
\tendL\mathcal{N} \left(\pmatrix{\mu_1 \cr\displaystyle  -\frac{1}{2}d_2^2
} , \pmatrix {d_1^2 & d_{11}\cr d_{11}& d_2^2
}\right) \qquad \mbox{as } \ny,
\end{equation}
then, under ${\mathrm{P}}\n_{\theta, \sigma+ n^{-1/2}\tau, 0 ;g_1}$,
$\Delta_{f_1;3}^{(n)}(\kappa)\tendL\mathcal{N}(\mu_1 + d_{11},
d_1^2)$. In the present context,
Proposition~\ref{lanart1} yields that
%
\begin{equation}\label{Lambdasigma}
\log\biggl(\frac{\mathrm{d}{\mathrm{P}}\n_{\theta, \sigma+
n^{-1/2}\tau, 0 ;g_1}}
{\mathrm{d}{\mathrm{P}}\n_{\theta, \sigma, 0 ;g_1}}
\biggr) = \tau\Delta\n_{g_1;2}(\theta, \sigma, 0) -\frac{1}{2}\tau
^2\sigma^{-2}\bigl(\mathcal{J}(g_1)-1\bigr) + \mathrm{o}_{\mathrm{P}}(1)
\end{equation}
as $\ny$, under ${\mathrm{P}}_{\theta, \sigma, 0;g_1}^{(n)}$; hence
(still as $\ny$ under ${\mathrm{P}}_{\theta, \sigma, 0;g_1}^{(n)}$),
%
\begin{eqnarray} \label{jointfgsigma}
\pmatrix{
\Delta_{f_1;3}^{(n)}(\kappa)\cr
\Delta\n_{g_1;2}(\theta, \sigma, 0)
}
&=&
n^{-1/2} \sum_{i=1}^n \pmatrix{
\phi_{f_1}(Z_i) ( Z_i^2 - \kappa )
\cr
\sigma^{-1}\bigl(\phi_{g_1}(Z_i) Z_i - 1\bigr)
}\nonumber\\[-8pt]\\[-8pt]
&\tendL&\mathcal{N} \left(
\pmatrix{
\mathcal{J}_{g_1}(f_1 )
\cr
0
} , \pmatrix{
d_1^2 & d_{11}
\cr
d_{11}& d_2^2
}\right),\nonumber
\end{eqnarray}
with $d_2^2=\sigma^{-2}(\mathcal{J}(g_1)-1)$ and
$
d_{11}
= \tau\sigma^{-1} \int_{-\infty}^\infty\phi_{f_1}(z) ( z^2 -
\kappa )(\phi_{g_1}(z) z - 1)g_1(z)\,\mathrm{d}z=0
$
(as the integral of a skew-symmetric function). We conclude that (\ref
{d11}) holds with $d_{11}=0$.

LeCam's third lemma therefore shows that the effect on the asymptotic
distribution of $\Delta_{f_1;3}^{(n)}(\kappa)$ of a $\mathrm{O}(n^{-1/2})$
perturbation of $\sigma$ is asymptotically nil; the asymptotic
linearity result of Proposition \ref{lin1art1} and a classical
argument on asymptotically discrete estimators (see,
e.g., Lemma 4.4 in \cite{K1987}) allow for extending this
conclusion to the stochastic $\mathrm{O}_{\mathrm{P}}(n^{-1/2})$ perturbations
induced by substituting a duly discretized root-$n$-consistent
estimator $\hat\sigma\n_\#$ for $\sigma$. Such a substitution
consequently does not affect the asymptotic behavior of $\Delta
_{f_1;3}^{(n)}(\kappa)$.

For $f_1\in\mathcal{F}_1$ and $g_1\in\mathcal{F}_{f_1}:=\{g_1\in
\mathcal{F}_1
\dvtx  \mathcal{K}_{g_1}(f_1) < \infty
\} $ (due to strong unimodality, $\mathcal{K}_{g_1}(f_1)<\infty$ also
implies $\mathcal{I}_{g_1}(f_1)<\infty$ and $\mathcal
{J}_{g_1}(f_1)<\infty$), let
%
\begin{equation}\label{gamman}\gamma\n(f_1) = \gamma\n(f_1, \theta
, \sigma):=\mathcal{K}\n(f_1) - 2\kappa(f_1) \mathcal{J}\n (f_1 )
+ \kappa^2(f_1) \mathcal{I}\n (f_1),
\end{equation}
where
%
\begin{eqnarray}\label{infog(f)1hat}\mathcal{I}\n(f_1) &=& \mathcal
{I}\n(f_1, \theta, \sigma):= n^{-1}\sum_{i=1}^n\phi^2_{f_1}(Z_i(
{\theta} , {\sigma} )),
\\ \label{infog(f)2hat}
\mathcal{J}\n(f_1 )&=& \mathcal{J}\n(f_1, \theta, \sigma
):=n^{-1}\sum_{i=1}^nZ_i^2( {\theta} , {\sigma} )\phi^2_{f_1}(Z_i(
{\theta} , {\sigma} ) )
\end{eqnarray}
and
%
\begin{equation}\label{infog(f)3hat}
\mathcal{K}\n(f_1)= \mathcal{K}\n(f_1, \theta, \sigma):=
n^{-1}\sum_{i=1}^n Z_i^4( {\theta} , {\sigma} ) \phi_{f_1}^2 (Z_i(
{\theta} , {\sigma} ))
\end{equation}
under ${\mathrm{P}}_{\theta, \sigma, 0;g_1}^{(n)}$ are consistent
estimates of $\mathcal{I}_{g_1} (f_1)$, $\mathcal{J}_{g_1} (f_1)$ and
$\mathcal{K}_{g_1} (f_1)$, respectively. Now in practice, $ \mathcal
{I}\n (f_1)$, $ \mathcal{J}\n (f_1 )$ and $\mathcal{K}\n(f_1) $,
hence $\gamma\n(f_1),$ cannot be computed from the observations and
$Z_i( {\theta} , {\hat\sigma}_\# )$ is to be substituted for $Z_i(
{\theta} , {\sigma} )$ in (\ref{infog(f)1hat})--(\ref
{infog(f)3hat}), 
yielding 
$ \gamma\n(f_1, {\theta} , {\hat\sigma}_\# )$. This substitution
in general requires a slight reinforcement of regularity assumptions.
Along the same lines as above (LeCam's third lemma and asymptotic
linearity), we easily obtain that $ \gamma\n(f_1, {\theta} , {\hat
\sigma}_\# )-\gamma\n(f_1, {\theta} , {\sigma})$ is $\mathrm{o}_{\mathrm{P}}(1)$
under ${\mathrm{P}}_{\theta, \sigma, 0;g_1}^{(n)}$ provided that the
asymptotic covariance of $\gamma\n(f_1, {\theta} , {\sigma})$ and
$\Delta_{g_1;2}^{(n)}$ is finite. A simple computation (and the strong
unimodality of $f_1$ and $g_1$) shows that a sufficient condition for
this is
%
\begin{equation}\label{casuffit} g_1\in\mathcal{F}_{f_1}^*:= \biggl\{
h_1\in\mathcal{F}_{f_1} \dvtx \int_{-\infty}^\infty z^5 \phi
_{f_1}^2(z) {\phi}_{h_1}(z) h_1(z) \,\mathrm{d}z <\infty \biggr\}.
\end{equation}
%
%
Defining the test statistic
%
\begin{equation}\label{stat2aart1}
\hat{T}^{(n)}_{f_1} ( {\theta}, {\sigma}):= \frac{1}{\sqrt{n\gamma
\n(f_1,\theta, \sigma)}}\sum_{i=1}^{n} \phi_{f_1}( Z_i( {\theta},
{\sigma}) ) \bigl( Z_i^2( {\theta} , {\sigma} ) - \kappa(f_1) \bigr)
\end{equation}
and the cross-information quantities
\begin{eqnarray*}
\mathcal{I}_{g_1}(f_1, g_1) &:=& \int_{- \infty}^{+ \infty} \phi
_{f_1} (z) \phi_{g_1} (z) g_1(z) \,\mathrm{d}z
,\\
\mathcal{J}_{g_1}(f_1, g_1) &:=& \int_{- \infty}^{+ \infty} z^2 \phi
_{f_1} (z) \phi_{g_1} (z) g_1(z) \,\mathrm{d}z
\end{eqnarray*}
and
\[
\mathcal{K}_{g_1}(f_1, g_1) := \int_{- \infty}^{+ \infty} z^4 \phi
_{f_1} (z) \phi_{g_1} (z) g_1(z) \,\mathrm{d}z
\]
(which for $f_1\in\mathcal{F}_1$ and $g_1\in\mathcal{F}^*_{f_1}$
are finite because of Cauchy--Schwarz), we thus have the following result.

\begin{Lem}\label{test2art1}
Let $f_1\in\mathcal{F}_1$ and $g_1\in\mathcal{F}_{f_1}^* $. Then:
\begin{enumerate}[(ii)]
\item[(i)] $\hat{T}^{(n)}_{f_1}( {\theta}, {\hat\sigma}_\#)=\hat
{T}^{(n)}_{f_1}( {\theta}, {\sigma})+\mathrm{o}_{\mathrm{P}}(1)$ is asymptotically
normal, with mean zero under ${\mathrm{P}}\n_{\theta, \sigma, 0;g_1}
$, mean
%
\begin{equation}\label{badshift}
\tau\  \frac{\mathcal{K}_{g_1}(f_1, g_1)-\mathcal{J}_{g_1}(f_1,
g_1)(\kappa(f_1)+\kappa(g_1)) + \mathcal{I}_{g_1}(f_1, g_1)\kappa
(f_1)\kappa(g_1)}
{[\mathcal{K}_{g_1}(f_1)-2\mathcal{J}_{g_1}(f_1)\kappa(f_1) +
\mathcal{I}_{g_1}(f_1)\kappa^2(f_1) ]^{1/2}}
\end{equation}
under $
{\mathrm{P}}_{\theta,\sigma,n^{-1/2}\tau; g_1}^{(n)}
$ and variance one under both.
\item[(ii)] The sequence of tests rejecting the null hypothesis
$\mathcal{H}\n_{\theta}:=\bigcup_{g_1\in\mathcal
{F}_{f_1}^*}\mathcal{H}\n_{\theta;g_1}$ of symmetry with respect to
specified $\theta$ whenever
$\hat{T}^{(n)}_{f_1}( {\theta}, \hat{\sigma}_\#)$ exceeds the
$(1-\alpha)$ standard normal quantile $z_\alpha$ is locally uniformly
asymptotically most powerful at asymptotic level $\alpha$ for
$\mathcal{H}\n_{\theta}$ against $\bigcup_{\xi>0} \bigcup_{
\sigma\in\R^+_0}\{{\mathrm{P}}\n_{\theta, \sigma, \xi;f_1}\}$.
\end{enumerate}
\end{Lem}
%


The tests based on $\hat{T}^{(n)}_{f_1}( {\theta}, \hat{\sigma}_\#
)$ enjoy all the validity (under $\mathcal{H}\n_{\theta}$) and
optimality (against $\bigcup_{\xi>0} \bigcup_{ \sigma\in\R^+_0}\{
{\mathrm{P}}\n_{\theta, \sigma, \xi;f_1}\}$) properties one can expect.
However, a closer look reveals that they are quite unsatisfactory on
one count: under $g_1\neq f_1$, their behavior strongly depends on the
arbitrary choice of the concept of scale (here, the median of absolute
deviations).

Consider, for example, the Gaussian version of (\ref{stat2aart1}),
which takes the form\vspace*{-3pt}
\[
\hat{T}^{(n)}_{\phi_1} ( {\theta}, {\sigma})= \frac{1}{\sqrt
{n\gamma\n(\phi_1)}}\sum_{i=1}^{n} \bigl( a Z_i^3( {\theta} , {\sigma}
) - 3 Z_i( {\theta} , {\sigma} ) \bigr),\vspace*{-3pt}
\]
where\vspace*{-3pt}
\[
\gamma\n(\phi_1) = \gamma\n(\phi_1 , \theta, \sigma)= a^2\sigma
^{-6}m_6^{(n)}(\theta) - 6 a \sigma^{-4}m_4^{(n)}(\theta) + 9\sigma
^{-2}m_2^{(n)}(\theta).\vspace*{-3pt}
\]
The test based on $\hat{T}^{(n)}_{\phi_1}( {\theta}, \hat{\sigma
})$ (here again, Slutsky's lemma allows for not
discretizing $ \hat{\sigma}$) is a pseudo-Gaussian test, hence optimal under Gaussian
assumptions; the asymptotic shift (\ref{badshift}) is $\tau\sqrt
{6/a}$ under $
{\mathrm{P}}_{\theta,\sigma,n^{-1/2}\tau; \phi_1}^{(n)} $, and
\[
\tau\bigl[5a\mu_4(g_1)-\bigl(9+3a\kappa(g_1)\bigr)\mu_2(g_1)+3\kappa(g_1)\bigr]
{[a^2 \mu_6(g_1)-6a\mu_4(g_1)+9\mu_2(g_1) ]^{-1/2}},
\]
where $\mu_k(g_1 ):=\int_{-\infty}^\infty z^k g_1(z) \,\mathrm{d}z$,
under $ {\mathrm{P}}_{\theta,\sigma,n^{-1/2}\tau; g_1}^{(n)}$. This
asymptotic shift strongly depends on $a$, hence on our (arbitrary)
choice of a scale parameter. Setting to one the standard deviation
instead of the median of absolute deviations would significantly modify
the local behaviour of $\hat{T}^{(n)}_{f_1}( {\theta}, \hat{\sigma
}_\#)$ as soon as $g_1\neq f_1$.
This does not affect optimality properties (which hold under $f_1$),
but is highly undesirable. 

Now, that unpleasant feature of $\hat{T}^{(n)}_{\phi_1}$ is entirely
due to the choice of $\kappa=\kappa(f_1)$ as a (nonrandom) centering
in (\ref{stat2aart1}). That choice was entirely motivated by
asymptotic orthogonality considerations under ${\mathrm{P}}_{\theta, \sigma
, 0;f_1}^{(n)}$, and does not affect the validity of the test. It
follows that replacing $\kappa(f_1)$ with any data-dependent sequence
$\kappa^{(n)}$ such that $\kappa\n-\kappa(f_1)=\mathrm{o}_{\mathrm{P}}(1)$ under
${\mathrm{P}}_{\theta, \sigma, 0;f_1}^{(n)}$ asymptotically has no impact
on $\hat{T}^{(n)}_{f_1}( {\theta}, {\sigma})$ under ${\mathrm{P}}_{\theta
, \sigma, 0;f_1}^{(n)}$. Let us show that this sequence $\kappa
^{(n)}$ can be chosen in order to cancel the unpleasant dependence of
the test statistic on the definition of scale.

Provided that\vspace*{-3pt}
\[
f_1\in\mathcal{F}_1^{\circ} :=\{h_1\in\mathcal{F}_1 \dvtx
z\mapsto\phi_{h_1}(z) \mbox{ is differentiable, with derivative }
\dot{\phi}_{h_1} \},\vspace*{-3pt}
\]
integration by parts yields\vspace*{-3pt}
\[
\mathcal{I}_{g_1}(f_1, g_1 )= \int_{-\infty}^{\infty} \dot{\phi
}_{f_1}(z) g_1(z)\,\mathrm{d}z\vspace*{-3pt}
\]
and\vspace*{-3pt}
\[
\mathcal{J}_{g_1}(f_1, g_1 )=2\int_{-\infty}^{\infty} z\phi
_{f_1}(z) g_1(z)\,\mathrm{d}z + \int_{-\infty}^{\infty} z^2\dot{\phi
}_{f_1}(z) g_1(z)\,\mathrm{d}z.\vspace*{-3pt}
\]
Therefore, $\mathcal{I}_{g_1}(f_1, g_1 )$, $\mathcal{J}_{g_1}(f_1,
g_1 )$ and $\kappa_{g_1}(f_1, g_1):= \mathcal{J}_{g_1}(f_1,g_1
)/\mathcal{I}_{g_1}(f_1,g_1)
$
under ${\mathrm{P}}_{\theta, \sigma, 0;g_1}^{(n)}$
are consistently estimated by\vspace*{-3pt}
\begin{eqnarray*}\label{In1}\mathcal{I}^{(n){\circ} }(f_1 ) =
\mathcal{I}^{(n){\circ} }(f_1, \theta,\sigma )
& := &
\frac{1}{n}\sum_{i=1}^n\dot{\phi}_{f_1}(Z_i(\theta, \sigma)) , \\
\mathcal{J}^{(n){\circ} }(f_1 ) =\mathcal{J}^{(n){\circ} }(f_1, \theta,\sigma )
& := &
\frac{2}{n}\sum_{i=1}^nZ_i(\theta, \sigma)\phi_{f_1}(Z_i(\theta,
\sigma) )
+\frac{1}{n}\sum_{i=1}^nZ_i^2(\theta, \sigma)\dot{\phi
}_{f_1}(Z_i(\theta, \sigma))\vspace*{-3pt}
\label{Jn1}
\end{eqnarray*}
and\vspace*{-2pt}
%
\begin{equation}\label{kappa*}\kappa^{(n){\circ} }(f_1) =
\kappa^{(n){\circ} }(f_1, \theta,\sigma ) :=\mathcal
{J}^{(n){\circ} }(f_1 )/ \mathcal{I}^{(n){\circ} }(f_1 ),\vspace*{-2pt}
\end{equation}
respectively.

Clearly, $\kappa^{(n){\circ} }(f_1)$ satisfies the requirement
that $\kappa^{(n){\circ} }(f_1)-\kappa(f_1)=\mathrm{o}_{\mathrm{P}}(1)$
under ${\mathrm{P}}_{\theta, \sigma, 0;f_1}^{(n)}$. In practice, however,
$ \kappa^{(n){\circ} }(f_1, \theta,\sigma )$ cannot be
computed from the observations, and $ \kappa^{(n){\circ}
}(f_1, \theta,{\hat\sigma}_\# )$, where $Z_i( {\theta} , {\hat
\sigma}_\# )$ has been substituted for $Z_i( {\theta} , {\sigma} )$,
is to be used instead. As in the evaluation of $\gamma\n(f_1)$
above (see (\ref{gamman})),
this substitution requires mild additional regularity conditions.
LeCam's third lemma then applies exactly along the same lines, implying that\vspace*{-2pt}
\[
\kappa^{(n){\circ} }(f_1, \theta,{\hat\sigma}_\# )-\kappa
^{(n){\circ} } (f_1, \theta,{ \sigma} ) = \mathrm{o}_{\mathrm{P}}(1)\vspace*{-2pt}
\]
under ${\mathrm{P}}_{\theta, \sigma, 0;g_1}^{(n)}$ as soon as the
asymptotic covariances of $\mathcal{I}^{(n){\circ} }(f_1 )$\vspace*{1pt}
and $\mathcal{J}^{(n){\circ} }(f_1 )$ with $\Delta
_{g_1;2}^{(n)}$ are finite. A simple computation (and the strong
unimodality of $f_1$ and $g_1$) shows that a sufficient conditions for
this is\vspace*{-2pt}
%
\begin{eqnarray} \label{casuffit2}
g_1\in\mathcal{F}_{f_1}^{\circ}&:=& \biggl\{ h_1\in\mathcal
{F}_{f_1}^{*} \dvtx \int_{-\infty}^\infty z^3 \dot{\phi}_{f_1}(z)
{\phi}_{h_1}(z) h_1(z) \,\mathrm{d}z <\infty \nonumber\\[-9pt]\\[-9pt]
&&\hphantom{\biggl\{}{}\mbox{and } \int_{-\infty
}^\infty z \dot{\phi}_{f_1}(z) {\phi}_{h_1}(z) h_1(z) \,\mathrm{d}z <\infty
\biggr\}\nonumber\vspace*{-2pt}
\end{eqnarray}
(no redundancy, since $ \dot{\phi}_{f_1}$ is not necessarily monotone).

Emphasize the dependence of $\Delta_{f_1;3}^{(n)}(\kappa)$ on $\theta
$ and $\sigma$ by writing\vspace*{1pt} $\Delta_{f_1;3}^{(n)}(\kappa, \theta,
\sigma)$: it follows from Lemma \ref{subst} in the \hyperref[appendixart1]{Appendix} that, for
$f_1\in\mathcal{F}_1^{\circ}$ and $g_1\in\mathcal
{F}_{f_1}^{\circ} $, the difference between $\Delta
_{f_1;3}^{(n)}(\kappa^{(n){\circ} } (f_1, {\theta},\hat
{\sigma}_\# ), {\theta}, \hat{\sigma}_\# )$ and $\Delta
_{f_1;3}^{(n)}(\kappa_{g_1}(f_1, g_1), \theta, {\sigma} )$ is
$\mathrm{o}_{\mathrm{P}}(1)$
under ${\mathrm{P}}_{\theta, \sigma, 0;g_1}^{(n)}$.
Letting (still for $f_1\in\mathcal{F}_1^{\circ}$)\vspace*{-2pt}
%
\begin{equation}\label{stat2bart1}
{T}_{f_1}^{(n){\circ} } ( {\theta}, {\sigma}) := \frac
{1}{\sqrt{n\gamma^{(n){\circ} }(f_1)}}\sum_{i=1}^{n} \phi
_{f_1}( Z_i( {\theta}, {\sigma}) ) \bigl( Z_i^2( {\theta} , {\sigma} ) -
\kappa^{(n){\circ} }(f_1) \bigr),\vspace*{-2pt}
\end{equation}
where\vspace*{-2pt}
\[
\gamma^{(n){\circ} }(f_1)=\gamma^{(n){\circ}
}(f_1,\theta, \sigma):=\mathcal{K}\n(f_1) - 2\kappa^{(n){\circ} } (f_1) \mathcal{J}\n (f_1 ) +\bigl( \kappa^{(n){\circ}
}(f_1) \bigr)^2 \mathcal{I}\n (f_1),\vspace*{-2pt}
\]
we thus have the following result.

\begin{Prop}\label{test2bart1}
Let $f_1\in\mathcal{F}_1^{\circ}$ and $g_1\in\mathcal
{F}_{f_1}^{\circ} $. Then:
\begin{enumerate}[(ii)]
\item[(i)] $ {T}^{(n){\circ} }_{f_1}( \theta, \hat{\sigma
}_\#)= {T}^{(n){\circ} }_{f_1}( {\theta}, {\sigma}) +\mathrm{o}_{\mathrm{P}}(1)$ is asymptotically normal, with mean zero under $ {\mathrm{P}}\n
_{\theta, \sigma, 0;g_1}$, mean\vspace*{-2pt}
%
\begin{equation}\label{noncentral}
\tau \frac{\mathcal
{K}_{g_1}(f_1, g_1)-\mathcal{J}_{g_1}(f_1, g_1)\kappa_{g_1}(f_1, g_1) }
{[\mathcal{K}_{g_1}(f_1)-2\mathcal{J}_{g_1}(f_1)\kappa_{g_1}(f_1,
g_1) + \mathcal{I}_{g_1}(f_1)\kappa^2 _{g_1}(f_1, g_1) ]^{1/2}}\vspace*{-2pt}
\end{equation}
under $ {\mathrm{P}}_{\theta,\sigma,n^{-1/2}\tau; g_1}^{(n)}$ and
variance one under both.
\item[(ii)] The sequence of tests rejecting the null hypothesis
$\mathcal{H}\n _\theta :=\bigcup_{g_1\in\mathcal{F}_{f_1}^{\circ}}\mathcal{H}\n_{\theta;g_1}$ of symmetry (with specified
location $\theta$, unspecified scale $\sigma$ and unspecified
standardized density $g_1\in\mathcal{F}_{f_1}^{\circ}$) whenever
$ {T}^{(n){\circ} }_{f_1}( \theta, \hat{\sigma}_\#) $
exceeds the $(1-\alpha)$ standard normal quantile $z_\alpha$ is
locally asymptotically most powerful at asymptotic level $\alpha$ for
$\mathcal{H}\n_\theta$ against $\bigcup_{\xi>0}\bigcup_{ \sigma
\in\R^+_0}\{{\mathrm{P}}\n_{\theta, \sigma, \xi;f_1}\}$.
\end{enumerate}
\end{Prop}

%

The advantage of the test statistic (\ref{stat2bart1}) compared to
(\ref{stat2aart1}) is that, irrespective of the underlying density
$g_1$, its behavior does not depend on the definition of the scale
parameter. The case of a Gaussian reference density ($f_1=\phi_1$),
however, is slightly different due to the particular form of the score
function $\phi_{f_1}$; see Section \ref{testspseudogauss}.

\subsubsection{Unspecified location}\label{322}
We now turn to the
case under which both $f_1$ and the location $\theta$ are unspecified.
Again, $\theta$ is to be replaced with some estimator, but additional
care has to be taken about the asymptotic impact of this substitution.
Still, from LeCam's third lemma, it follows that the impact, under
${\mathrm{P}}_{\theta, \sigma, 0;g_1}^{(n)}$, of an estimated $\theta$ on
$\Delta_{f_1;3}^{(n)}(\kappa)$ can be obtained from the asymptotic
behavior of
\[
\pmatrix{
\Delta_{f_1;3}^{(n)}(\kappa)\cr
\Delta\n_{g_1;1}(\theta, \sigma, 0)
}
=
n^{-1/2} \sum_{i=1}^n \pmatrix{
\phi_{f_1}(Z_i) ( Z_i^2 - \kappa )\cr
\sigma^{-1}\phi_{g_1}(Z_i)
},
\]
which is asymptotically normal with asymptotic covariance matrix
\[ \pmatrix{
\gamma_{g_1}^\kappa(f_1)&\delta_{g_1}^\kappa(f_1,g_1) \cr
\delta_{g_1}^\kappa(f_1,g_1)& \sigma^{-2}\mathcal{I}(g_1 )
},
\]
where
$\delta_{g_1}^\kappa(f_1,g_1):= \sigma^{-1}( \mathcal
{J}_{g_1}(f_1,g_1) -\kappa \mathcal{I}_{g_1}(f_1,g_1))
.$
Clearly, this covariance $\delta_{g_1}^\kappa(f_1,g_1)$ vanishes iff
$\kappa=\kappa_{g_1}(f_1, g_1),$ which, for $g_1=f_1$, coincides with
$\kappa(f_1)$.

Assuming that an estimate $\kappa\n(f_1)$ such that $\kappa\n
(f_1)-\kappa_{g_1}(f_1, g_1)=\mathrm{o}_{\mathrm{P}}(1)$ under ${\mathrm{P}}_{\theta,
\sigma, 0;g_1}^{(n)}$ exists, $\Delta_{f_1;3}^{(n)}(\kappa\n(f_1)
)$ is asymptotically equivalent to $\Delta_{f_1;3}^{(n)}(\kappa(f_1)
)$ under~${\mathrm{P}}_{\theta, \sigma, 0;f_1}^{(n)}$, and asymptotically
uncorrelated with $\Delta\n_{g_1;1}(\theta, \sigma, 0)$ and $\Delta
\n_{g_1;2}(\theta, \sigma, 0)$ (hence, asymptotically insensitive
(in probability) to root-$n$ perturbations of both $\theta$ and
$\sigma$)
under ${\mathrm{P}}_{\theta, \sigma, 0;g_1}^{(n)}$. It follows from
Section \ref{321} that $\kappa^{(n){\circ}}(f_1,\theta,
\sigma)$ defined in (\ref{kappa*}) is such an estimator. The same
reasoning as in Section \ref{321} implies that this still holds when
substituting, in $\Delta_{f_1;3}^{(n)}(\kappa)$, any estimators $\hat
\theta_\#$ and $\hat\sigma_\#$ satisfying (C1) and (C2) for $\theta
$ and~$\sigma$.
Finally, Lemma \ref{subst} in the \hyperref[appendixart1]{Appendix} ensures that $\Delta
_{f_1;3}^{(n)}(\kappa^{(n){\circ} } (f_1, \hat{\theta}_\#
,\hat{\sigma}_\# ), \hat{\theta}_\#, \hat{\sigma}_\# )$ can be
substituted for $\Delta_{f_1;3}^{(n)}(\kappa_{g_1}(f_1, g_1), \theta
, {\sigma} )$.
We thus have shown the following result.

\begin{Prop}\label{test3art1}
Let $f_1\in\mathcal{F}_1^{\circ}$ and $g_1\in\mathcal
{F}_{f_1}^{\circ}$. Then:
\begin{enumerate}[(ii)]
\item[(i)] $ {T}^{(n){\circ}}_{f_1}( \hat{\theta}_\# , \hat
{\sigma}_\# )={T}^{(n){\circ}}_{f_1}( {\theta} , \hat{\sigma
}_\# )
+\mathrm{o}_{\mathrm{P}}(1) ={T}^{(n){\circ}}_{f_1}( {\theta}, {\sigma
})+\mathrm{o}_{\mathrm{P}}(1)$ is asymptotically normal, with mean zero under $ {\mathrm{P}}\n_{\theta, \sigma, 0;g_1}$, mean (\ref{noncentral})
under $ {\mathrm{P}}_{\theta,\sigma,n^{-1/2}\tau; g_1}^{(n)}$ and
variance one under both.
\item[(ii)] The sequence of tests rejecting the null hypothesis
$\mathcal{H}\n:=\bigcup_{g_1\in\mathcal{F}_{f_1}^{\circ}}\bigcup_{\theta\in\R}\mathcal{H}\n_{\theta;g_1}$ of symmetry
(with unspecified location $\theta$, unspecified scale $\sigma$ and
unspecified standardized density $g_1$) whenever
$ {T}^{(n){\circ}}_{f_1}( \hat{\theta}_\# , \hat{\sigma}_\#
)$ exceeds the $(1-\alpha)$ standard normal quantile $z_\alpha$ is
locally asymptotically most powerful at asymptotic level $\alpha$ for
$\mathcal{H}\n$ against $\bigcup_{\xi>0}\bigcup_{\theta\in\R
}\bigcup_{ \sigma\in\R^+_0}\{{\mathrm{P}}\n_{\theta, \sigma, \xi
;f_1}\}$.
\end{enumerate}
\end{Prop}

This test is based on the same test statistic $ {T}_{f_1}^{(n){\circ}} $ as the specified-location test of
Proposition~\ref{test2bart1}, except that the (here unspecified) location $\theta$ is
replaced by an estimator~$ \hat{\theta}_\#$. The local powers of the
two tests coincide: asymptotically, again, there is no loss of
efficiency due to the non-specification of $\theta$.

\subsection{Pseudo-Gaussian tests }\label{testspseudogauss}
Particularizing the reference density $f_1$ as the standard normal one
$\phi_1$ in the tests of Sections~\ref{321} and \ref{322} in
principle yields \textit{pseudo-Gaussian} tests, based on the test
statistics $T_{\phi_1}^{(n){\circ}} ( {\theta})$ or $T_{\phi
_1}^{(n){\circ}} ( {\hat\theta})$. Due to the particular form
of the Gaussian score function, however, the Gaussian statistic can be
given a much simpler form.
Indeed, $\mathcal{I}_{g_1}(\phi_1, g_1)=\mathcal{I}(\phi_1) = a$
does not depend on $g_1$, and needs not be estimated,
while $\mathcal{J}_{g_1}(\phi_1, g_1)=\mathcal{J}(\phi_1) = 3a\mu
_2(g_1)$, so that $\kappa_{g_1}(f_1, g_1)$ is consistently estimated
by $3m_2^{(n)}(\theta)/\sigma^2$. This, after elementary computation,
yields the test statistic
%
\begin{equation}\label{statGbart1}
{T}^{(n)\dagger}(\theta) :=\frac{1}{\sqrt{n \gamma^{(n)\dagger}
}}\sum_{i=1}^{n} (X_i-\theta) \bigl( (X_i-\theta)^2 - 3 m_2^{(n)}(\theta
) \bigr),
\end{equation}
where $ \gamma^{(n)\dagger} := \gamma^{(n)\dagger}(\theta)
:=m_6^{(n)}(\theta) - 6m_2^{(n)}(\theta) m_4^{(n)}(\theta) +9(
m_2^{(n)}(\theta) )^3 $.
For this test statistic ${T}^{(n)\dagger}(\theta)$, the asymptotic
shift (\ref{noncentral}) under $ {\mathrm{P}}_{\theta,\sigma,n^{-1/2}\tau
; g_1}^{(n)}$ now takes the form
\[
\tau[5\mu_4(g_1) -9\mu_2^2(g_1) ]
{[\mu_6(g_1) -6\mu_2(g_1) \mu_4(g_1) +9\mu_2^3(g_1) ]^{-1/2}} .
\]
This shift does not depend on $a$ anymore, and still reduces to\vspace*{-2pt} $\tau
\sqrt{6/a}$ under $ {\mathrm{P}}_{\theta,\sigma,n^{-1/2}\tau; \phi
_1}^{(n)}$ (the same value as for $\hat{T}^{(n)}_{\phi_1} ( {\theta
}, {\sigma})$, which confirms that optimality under Gaussian densities
has been preserved); nor does it depend on the scale.

The tests based on the asymptotically standard normal null distribution
of ${T}^{(n)\dagger}$ are optimal under Gaussian assumptions, but
remain valid when those assumptions are violated.
Again, a simple Slutsky argument allows for replacing $\theta$ (if
unspecified) with any consistent estimator $\hat{{\theta}}$ without
going through discretization; moreover, (\ref{statGbart1}) does not
depend on $\sigma$.
The tests based on ${T}^{(n)\dagger} ( {\theta})$ and
${T}^{(n)\dagger}( \bar{X}\n)$ both are closely related to the
traditional test of symmetry based on $b_{1}^{(n)}$. More precisely,
under any $ {\mathrm{P}}_{\theta,\sigma,0 ; g_1}^{(n)} $, $g_1\in
(\mathcal{F}_{\phi_1}^{\circ} =) \mathcal{F}_{\phi_1} $
(note that the assumption $g_1\in\mathcal{F}_{\phi_1}^{\circ} = \mathcal{F}_{\phi_1} $ implies that $g_1$ has finite moments of
order six),
%
\begin{equation}\label{S=T}
{T}^{(n)\dagger}( {\theta})={T}^{(n)\dagger} \bigl( \bar{X}\n\bigr) +\mathrm{o}_{\mathrm{P}}( 1 ) =
S\n_2 +\mathrm{o}_{\mathrm{P}}( 1 ),
\end{equation}
where $S\n_2$ is the empirically standardized form (\ref{classical})
of $b_{1}^{(n)}$ (see (\ref{classical})).

Summing up, we thus have the following result.

\begin{Prop}\label{prop3.4} 
Let $g_1\in\mathcal{F}_{\phi_1}$ and $\hat{\theta}=\theta+
\mathrm{O}_{\mathrm{P}}(n^{-1/2})$; recall that $\mu_k(g_1 ) :=\int_{-\infty
}^\infty z^k g_1(z) \,\mathrm{d}z$ stands for $ g_1$'s moment of order $k$.
Then:
\begin{enumerate}[(iii)]
\item[(i)] ${T}^{(n)\dagger}( \hat{\theta})={T}^{(n)\dagger} (
{\theta})+\mathrm{o}_{\mathrm{P}}(1)$ is asymptotically normal, with mean zero under
$ {\mathrm{P}}\n_{\theta, \sigma, 0;g_1}$, mean
\[
\tau[5\mu_4(g_1 )
-9\mu_2^2(g_1 ) ]/
[\mu_6(g_1 ) - 6\mu_2(g_1 ) \mu_4(g_1 ) + 9\mu_2^3(g_1 ) ]^{1/2}
\]
under $ {\mathrm{P}}_{\theta,\sigma,n^{-1/2}\tau; g_1}^{(n)}$ and
variance one under both.
\item[(ii)] The sequence of tests rejecting the null hypothesis of
symmetry (with specified location $\theta$) $\mathcal{H}\n_{\theta}
:=\bigcup_{g_1\in\mathcal{F}_{\phi_1}}\mathcal{H}\n_{\theta
;g_1}$ whenever
$ {T}^{(n)\dagger} ( {\theta})$ exceeds the $(1-\alpha)$ standard
normal quantile $z_\alpha$ is locally asymptotically most powerful at
asymptotic level $\alpha$ against $\bigcup_{\xi>0}\bigcup_{ \sigma
\in\R^+_0}\{{\mathrm{P}}\n_{\theta, \sigma, \xi;\phi_1}\}$.
\item[(iii)] The sequence of tests rejecting the null hypothesis of
symmetry (with unspecified location) $\mathcal{H}\n:=\bigcup_{g_1\in
\mathcal{F}_{\phi_1}}\bigcup_{\theta\in\R}\mathcal{H}\n_{\theta
;g_1}$ whenever
${T}^{(n)\dagger} ( \hat{\theta} )$ exceeds the $(1-\alpha)$
standard normal quantile $z_\alpha$ is locally asymptotically most
powerful at asymptotic level $\alpha$ against $\bigcup_{\xi
>0}\bigcup_{\theta\in\R}\bigcup_{ \sigma\in\R^+_0}\{{\mathrm{P}}\n
_{\theta, \sigma, \xi;\phi_1}\}$.
\end{enumerate}
\end{Prop}
%
%

For the sake of completeness, we also provide (with the same notation)
the following result on the asymptotic behavior of the (suboptimal)
test based on $m\n_3(\theta)$. Details are left to the reader.

\begin{Prop}\label{testGbart2}
Let $g_1\in\mathcal{F}_{\phi_1}$. Then,
$S\n_1 := n^{1/2}m\n_3(\theta)/(m\n_6(\theta))^{1/2}$ is
asymptotically normal, with mean zero under $ {\mathrm{P}}\n_{\theta,
\sigma, 0;g_1}$, mean $\tau[5\mu_4(g_1 ) -3\kappa(g_1) \mu_2(g_1 ) ]/
\mu_6^{1/2}(g_1 ) $
under $ {\mathrm{P}}_{\theta,\sigma,n^{-1/2}\tau; g_1}^{(n)}$ and
variance one under both.
\end{Prop}

Under Gaussian densities ($g_1 = \phi_1$), the asymptotic shifts of
$T_{\phi_1}^{(n){\circ}} ( {\theta})$ (Proposition~\ref{test3art1}(i))
and $S\n_1$ (Proposition~\ref{testGbart2}) are
$16\tau/\sqrt{6}
$ and $16\tau/\sqrt{15}$, respectively; the asymptotic relative
efficiency of ${T}^{(n)\dagger} ( {\theta})$ with respect to $S\n_1$
is thus as high as 2.5 in the vicinity of Gaussian densities. This,
which is not a small difference, confirms the suboptimality of $m\n
_3(\theta)$-based tests.

\subsection{Laplace tests}\label{testsLaplace}

Replacing the Gaussian reference density $\phi_1$ with the
double-exponential one~$f_\mathcal{L}$, we similarly obtain the \textit{Laplace tests}.
The assumption that $f_1\in\mathcal{F}_1^{\circ}$ unfortunately rules out~$f_\mathcal{L}$,
since $\phi_{f_\mathcal
{L}}(z)=\mbox{sign}(z)/d$ is not differentiable, so that the
\mbox{construction of $\kappa^{(n){\circ} }(f_1)$} in~(\ref{kappa*})
does not apply for $f_1= f_\mathcal{L}$. Now, a direct construction is
possible: $\mathcal{I}_{g_1}(f_\mathcal{L}, g_1)$ indeed reduces to
$2g_1(0)/d$ -- which is consistently estimated by $\mathcal
{I}^{(n){\circ} }(f_\mathcal{L} ):=2\hat{g}_1(0)/d$ (where
$\hat{g_1}$, e.g., is some kernel estimator of $g_1$).
Similarly, $\mathcal{J}_{g_1}(f_\mathcal{L}, g_1)$ reduces to
$({2}/{d})\int_{-\infty}^{\infty} |z| g_1(z) \,\mathrm{d}z$ -- which is
consistently estimated by $\mathcal{J}^{(n){\circ}
}(f_\mathcal{L} ):=({2}/{nd})\sum_{i=1}^n
|Z_i(\theta,\break \hat\sigma_\#)|$; the scaling constant $d$ is easily
computed, yielding $d=1/(\log2)\approx1.44$.
Then,
\[
\kappa^{(n){\circ} }(f_\mathcal{L}):=\mathcal{J}^{(n){\circ} }(f_\mathcal{L} )/\mathcal{I}^{(n){\circ}
}(f_\mathcal{L} )= \frac{1}{n \hat{g}_1(0)}\sum_{i=1}^n
|Z_i(\theta, \hat\sigma_\#)| 
\]
is such that $\kappa^{(n){\circ} }(f_\mathcal{L})-\kappa
(f_\mathcal{L})=\mathrm{o}_{\mathrm{P}}(1)$ under ${\mathrm{P}}_{\theta, \sigma,
0;f_\mathcal{L}}^{(n)}$, as required.

The Laplace tests are based on $T_\mathcal{L}^{(n){\circ}} (
{\theta})$ (specified $\theta$) or $T_\mathcal{L}^{(n){\circ}} ( {\hat\theta})$ (unspecified $\theta$), where
%
\begin{equation}\label{TLaplace}T_\mathcal{L}^{(n){\circ}} (
{\theta}):=
\frac{1}{\sqrt{n \gamma^{(n){\circ}} (f_\mathcal{L})}} \sum
_{i=1}^n \operatorname{sign}(Z_i(\theta, \hat\sigma_\#)) \bigl( (Z_i(\theta, \hat
\sigma_\#))^2-\kappa^{(n){\circ} }(f_\mathcal{L}) \bigr)
\end{equation}
with
\[
\gamma^{(n){\circ}} (f_\mathcal{L}) = \frac{m\n_4}{{\hat
\sigma}_\#^4}-2\frac{m\n_2}{n{\hat\sigma}_\#^2\hat{g}_1(0)}\sum
_{i=1}^n \vert Z_i(\theta, \hat\sigma_\#) \vert+ \Biggl(\frac{1}{n\hat
{g}_1(0)}\sum_{i=1}^n \vert Z_i(\theta, \hat\sigma_\#) \vert \Biggr)^2.
\]
%
These tests share with the Gaussian Fechner test (see \cite{CHP2008}) the use of the score function $z\mapsto\operatorname{sign}(z) z^2$. The
orthogonalization, however, differs, since the Fechner and Edgeworth
families the tests were built on are different.
The following proposition summarizes their properties; details are left
to the reader.

\begin{Prop}\label{testGbart1}\label{prop3.6}
Let $g_1\in\mathcal{F}_{f_\mathcal{L}}$, $\hat{\theta}=\theta+
\mathrm{O}_{\mathrm{P}}(n^{-1/2})$ and denote by $\mu_{|k|}(g_1 ) :=\int_{-\infty
}^{\infty} |z|^k g_1(z) \,\mathrm{d}z$ the absolute moment of order $k$ of $g_1$. Then,
\begin{enumerate}[(iii)]
\item[(i)] $ T_\mathcal{L}^{(n){\circ}}( \hat{\theta
})=T_\mathcal{L}^{(n){\circ}}( {\theta})+\mathrm{o}_{\mathrm{P}}(1)$ is
asymptotically normal, with mean zero under $ {\mathrm{P}}\n_{\theta,
\sigma, 0;g_1}$, mean
\[
\frac{\tau[4\mu_{|3|}(g_1 ) -2\mu_{|1|}^2(g_1 ) /g_1(0)]}
{[\mu_4(g_1 ) -2\mu_2(g_1 ) \mu_{|1|}(g_1 ) /g_1(0)+\mu_{|1|}^2(g_1
) /(g_1(0))^2 ]^{1/2}}
\]
under $ {\mathrm{P}}_{\theta,\sigma,n^{-1/2}\tau; g_1}^{(n)}$ and
variance one under both.
\item[(ii)] The sequence of tests rejecting the null hypothesis of
symmetry (with specified location $\theta$) $\mathcal{H}\n_{\theta}
:=\bigcup_{g_1\in\mathcal{F}_{f_\mathcal{L}}}\mathcal{H}\n
_{\theta;g_1}$ whenever
$ T_\mathcal{L}^{(n){\circ}}(\theta)$ exceeds the $(1-\alpha
)$ standard normal quantile $z_\alpha$ is locally asymptotically most
powerful at asymptotic level $\alpha$ against $\bigcup_{\xi
>0}\bigcup_{ \sigma\in\R^+_0}\{{\mathrm{P}}\n_{\theta, \sigma, \xi
;f_\mathcal{L}}\}$.
\item[(iii)] The sequence of tests rejecting the null hypothesis of
symmetry (with unspecified location) $\mathcal{H}\n:=\bigcup_{g_1\in
\mathcal{F}_{f_\mathcal{L}}}\bigcup_{\theta\in\R}\mathcal{H}\n
_{\theta;g_1}$ whenever
$ T_\mathcal{L}^{(n){\circ}} ( \hat{\theta} )$ exceeds the
$(1-\alpha)$ standard normal quantile $z_\alpha$ is locally
asymptotically most powerful at asymptotic level $\alpha$ against
$\bigcup_{\xi>0}\bigcup_{\theta\in\R}\bigcup_{ \sigma\in\R
^+_0}\{{\mathrm{P}}\n_{\theta, \sigma, \xi;f_\mathcal{L}}\}$.
\end{enumerate}
\end{Prop}

Comparing the asymptotic shifts of the pseudo-Gaussian tests
and the Laplace ones
yields asymptotic relative efficiency values; the asymptotic efficiency
of tests based on~${T}^{(n)\dagger}$ with respect to those based on
$T_\mathcal{L}^{(n){\circ}}$ is 1.76 in the vicinity of
Gaussian densities and 0.7 in the vicinity of double exponential ones.
Finally, note that the empirical me\-dian~$X\n_{1/2}$ here provides a
much more sensible estimator of $\theta$ than the empirical mean~$\bar{X}\n$;
it has been used for $\hat\theta$ in the simulations of
$T_\mathcal{L}^{(n){\circ}}({\hat\theta})$ in Section \ref{simu}.

\section{Finite sample performances}\label{simu}

We performed a first simulation study on the basis of $N=5000$
independent samples of size $n=100$ from (\ref{modelart1}), with
normal and double-exponential densities $f_1$ and
skewness parameter values $\xi=0.1$ and $\xi=0.2$. Each of those
samples was subjected, at asymptotic level $\alpha=5\%$, to the
classical specified location test of skewness based on $m\n_3(\theta
)$ (i.e., on (\ref{classicaltheta})), the (optimal) pseudo-Gaussian
tests based on $b\n_1$ (i.e., on (\ref{classical})) and the
corresponding Laplace and logistic tests. For the sake of completeness,
the two \textit{triples tests} proposed by Randles et al. \cite{RFPW1980}, which
are based on the signs of $X_i+X_j-2X_k$, $1\leq i<j<k\leq n$, are also
included in this simulation study. Those tests, which are location-invariant, do not follow from any argument of group invariance and are
not distribution-free.

Rejection frequencies are reported in Table \ref{tablesimu1}.

\begin{table}
\caption{Rejection frequencies (out of $N=5000$ replications),
under various symmetric and skewed normal and double-exponential
distributions from the Edgeworth families (\protect\ref{modelart1}), with $\xi
=0, 0.1, 0.2$, of the classical tests of skewness, based on $m\n
_3(\theta)$ and $b\n_1$, the Gaussian, Laplace and logistic tests,
and the \textit{triples tests} $T_{R1}\n$ and $T_{R1}\n$ of \protect\cite{RFPW1980}}\label{tablesimu1}
\begin{tabular*}{\tablewidth}{@{\extracolsep{\fill}}lllllll@{}}
\hline
Test& \multicolumn{3}{l}{$\mathcal{SN}(\xi)$}
& \multicolumn{3}{l@{}}{$\mathcal{SL}(\xi)$}\\[-5pt]
&\multicolumn{3}{c}{\hrulefill}&\multicolumn{3}{c@{}}{\hrulefill}\\
 & \multicolumn{3}{l}{$\xi$} &
\multicolumn{3}{l@{}}{$\xi$} \\[-5pt]
&\multicolumn{3}{c}{\hrulefill}&\multicolumn{3}{c@{}}{\hrulefill}\\
&0&0.1&0.2&0&0.1&0.2\\
\hline
$m_3^{(n)}(\theta)$ &0.0372 & 0.1136 & 0.0996 & 0.0306 & 0.6938 &
0.8722 \\[1pt]
${T}^{(n)\dagger} ( {\theta})$ & 0.0434 & 0.7276 & 0.9958 & 0.0252 &
0.4596 & 0.6774\\[1pt]
$b\n_1$ &0.0416 & 0.6986 & 0.9746 & 0.0444 & 0.7458 & 0.8930 \\[1pt]
$T_\mathcal{L}^{(n){\circ}} ( {\theta})$& 0.0520 & 0.5424 &
0.9474 & 0.0406 & 0.9090 & 0.9998 \\[1pt]
$T_\mathcal{L}^{(n){\circ}} (X\n_{1/2})$ & 0.0280 & 0.4440 &
0.8360 & 0.0284 & 0.8838 & 0.9960 \\[1pt]
$T_{\mathrm{Log}}^{(n){\circ}} ( {\theta})$ &0.0492 & 0.7336 &
0.9954 & 0.0378 & 0.8516 & 0.9894 \\[1pt]
$T_{\mathrm{Log}}^{(n){\circ}} ( \bar{X}\n)$ & 0.0362 & 0.6626 &
0.9716 & 0.0384 & 0.8516 & 0.9880 \\[1pt]
$T\n_{R1}$ &0.0518 & 0.6786 & 0.9606 & 0.0576 & 0.9276 & 0.9986 \\[1pt]
$T\n_{R2}$ &0.0608 & 0.6992 & 0.9640 & 0.0650 & 0.9350 & 0.9988 \\
\hline
\end{tabular*}
\end{table}

Note that all tests considered here, except for Randles', are extremely
conservative, and in most cases hardly reach the nominal 5\% rejection
frequency under the null. Randles' tests, on the other hand,
significantly over-reject, which does not facilitate comparisons.
Despite that, the tests based on $b\n_1$, $T_\mathcal{L}^{(n){\circ}} ( X\n_{1/2})$ and $T_{\mathrm{Log}}^{(n){\circ}} ( \bar
{X}\n)$ exhibit excellent peformances, and largely outperform those
based on $m\n_3(\theta)$ (although the latter requires~%
$\theta$ to be known).


The Edgeworth families considered throughout this paper, however,
served as a theoretical guideline in the construction of our Edgeworth
testing procedures, and never were meant as an actual data generating
process. One could argue that analyzing performances under alternatives
of the Edgeworth type creates an unfair bias in favor of our methods.
Therefore, we also generated $N=5000$ independent samples of size
$n=100$ from the skew-normal $\mathcal{SN}(\lambda)$ and skew-$t$
$\mathcal{S}t(\nu, \lambda)$ densities (with $\nu=2$, $\nu=4$ and
$\nu=8$ degrees of freedom) defined by Azzalini and Capitanio \cite{AC2003},
for various values of their skewness coefficient $\lambda$ ($\lambda
=0$ implying symmetry); since the sign of $\lambda$ is not directly
related to that of $\xi$, we only performed two-sided tests.
That class of skewed densities was chosen in view of its increasing
popularity among practitioners.

None of the tests considered in this simulation example are optimal in
this Azzalini and Capitanio context.
Inspection of Table \ref{tablesimu2} nevertheless reveals that the
classical tests of skewness based on $m_3^{(n)}(\theta)$ and $b\n_1$
collapse under $t_2$ and $t_4$, which have infinite sixth-order
moments, and under the related $\mathcal{S}t(2, \lambda)$ and
$\mathcal{S}t(4, \lambda)$ densities. The same tests fail to achieve
the 5\% nominal level under the Student distribution with 8 degrees of
freedom (despite finite sixth-order moments) and show weak performance
under the $\mathcal{S}t(8, \lambda)$ density. Remark that the
suboptimality of the test based on $m_3^{(n)}(\theta)$, which, as a
consequence of Proposition \ref{test1art1}, may be considered as an artificial
consequence of the choice of skewed families of the Edgeworth type,
nevertheless also very neatly appears here. The triples tests behave
uniformly well; note, however, their tendency to over-rejection, in
particular under Student densities.

\begin{sidewaystable}
\tabcolsep=0pt
\tablewidth=\textheight
\tablewidth=\textwidth
\caption{Rejection frequencies (out of $N=5000$ replications),
under various symmetric and related skew-normal and skew-$t$
distributions \protect\cite{AC2003} $\mathcal{SN}(\lambda)$
and $\mathcal{S}t(\nu,\lambda)$ ($\nu=2, 4, 8$ and various $\lambda
$) of the classical tests of skewness, based on $m\n_3(\theta)$ and
$b\n_1$, the Gaussian, Laplace and logistic tests, and the \textit{triples tests} $T_{R1}\n$ and $T_{R1}\n$ of \cite{RFPW1980}}
\label{tablesimu2}
\begin{tabular*}{\tablewidth}{@{\extracolsep{\fill}}lllllllllllllllll@{}}
\hline
Test&\multicolumn{16}{l@{}}{$\lambda$}\\[-5pt]
&\multicolumn{16}{c@{}}{\hrulefill}\\
 & \multicolumn{4}{l}{$\mathcal{SN}(\lambda
)$} & \multicolumn{4}{l}{$\mathcal{S}t(2,\lambda)$}  & \multicolumn{4}{l}{$\mathcal
{S}t(4,\lambda)$} & \multicolumn{4}{l@{}}{$\mathcal{S}t(8,\lambda)$}\\[-5pt]
&\multicolumn{4}{c}{\hrulefill}&\multicolumn{4}{c}{\hrulefill}&\multicolumn{4}{c}{\hrulefill}&\multicolumn{4}{c@{}}{\hrulefill}\\
 &0&1&2&3&0&2&4&6&0&2&4&6&0&2&4&6\\
\hline
$m_3^{(n)}(\theta)$ &0.0476 & 0.0482 & 0.0952 & 0.1936 & 0.0072 &
0.0106 & 0.0126 & 0.0144 &0.0192 & 0.0190 & 0.0298 & 0.0436 & 0.0316 & 0.0608 & 0.1302 & 0.1754 \\[2pt]
${T}^{(n)\dagger} ( {\theta})$ & 0.0374 & 0.0634 & 0.2988 & 0.5942 &
0.0046 & 0.0118 & 0.0182 & 0.0268& 0.0144 & 0.0184 & 0.0586 & 0.1134 & 0.0260 & 0.1422 & 0.4078 & 0.5428 \\
$b\n_1$ & 0.0418 & 0.0616 & 0.3066 & 0.6130 & 0.0172 & 0.0232 & 0.0308
& 0.0396 &0.0252 & 0.0302 & 0.0754 & 0.1298 & 0.0322 & 0.1604 & 0.4186 & 0.5484 \\[1pt]
$T_\mathcal{L}^{(n){\circ}} ( {\theta})$& 0.0460 & 0.0690 &
0.2682 & 0.5406 & 0.0180 & 0.0414 & 0.0740 & 0.0988 & 0.0332 & 0.0406 &0.1304 & 0.2514 & 0.0456 & 0.2054 & 0.5640 & 0.6906 \\[2pt]
$T_\mathcal{L}^{(n){\circ}} (X\n_{1/2})$ & 0.0334 & 0.0472 &
0.2022 & 0.4736 & 0.0168 & 0.0288 & 0.0520 & 0.0678 & 0.0206 & 0.0302 & 0.0950 & 0.1798 & 0.0304 & 0.1518 & 0.4834 & 0.6260 \\
$T_{\mathrm{Log}}^{(n){\circ}} ( {\theta})$ &0.0468 & 0.0742 &
0.3542 & 0.7010 & 0.0154 & 0.0286 & 0.0496 & 0.0666 &0.0236 & 0.0318 & 0.1082 & 0.2086 &0.0342 & 0.2104 & 0.5846 & 0.7508 \\[3pt]
$T_{\mathrm{Log}}^{(n){\circ}} ( \bar{X}\n)$ & 0.0354 & 0.0568 &
0.2988 & 0.6426 & 0.0144 & 0.0256 & 0.0408 & 0.0492 & 0.0228 & 0.0276 & 0.0882 & 0.1694 & 0.0288 & 0.1712 & 0.5046 & 0.6696 \\[3pt]
$T\n_{R1}$ &0.0540 & 0.0778 & 0.3602 & 0.7082 & 0.0618 & 0.1032 &
0.1798 & 0.2312 & 0.0508 & 0.0636 & 0.1842 & 0.3444 & 0.0530 & 0.2592 & 0.6766 & 0.8336 \\[2pt]
$T\n_{R2}$ &0.0598 & 0.0886 & 0.3812 & 0.7258 & 0.0656 & 0.1098 &
0.1882 & 0.2424 & 0.0556 & 0.0688 & 0.1938 & 0.3582 & 0.0598 & 0.2740 & 0.6940 & 0.8422 \\
\hline
\end{tabular*}
\end{sidewaystable}

\section{Conclusions and perspectives}\label{conclu}

We have derived the optimal tests for testing the hypothesis of
symmetry within families of skewed densities mimicking, in the Gaussian
case, the type of local asymmetry observed in a central limit behaviour.
The resulting tests were obtained for specified and unspecified
locations, under specified or unspecified densities (satisfying
appropriate assumptions). Local powers and asymptotic relative
efficiencies are computed and finite-sample performances are
investigated in the context of classical skew-normal and skew-$t$ distributions.

Establishing the optimality properties of the traditional test based on
the Pearson--Fisher coefficient $b\n_1$ was one of the objectives of
this work, and we show that this test is indeed optimal in the vicinity
of Gaussian symmetry. Interestingly, its optimality holds for the
specified-location as well as for the unspecified-location hypothesis
of symmetry and is preserved if centering, in the computation of the
test statistic, is based on robust root-$n$-consistent estimators of
location rather than on the sample mean ${\bar X}\n$.

Table \ref{summary} provides a summary of the various tests described
throughout Section 3, along with their validity and optimality properties.

\begin{sidewaystable}
\tabcolsep=0pt
\tablewidth=\textheight
\tablewidth=\textwidth
\caption{A summary of the various test statistics considered
throughout the paper, with reference to their definitions, the testing
problem (specified or unspecified location -- scale throughout remains
unspecified) addressed, the standardized densities $g_1$ for which
they are valid, the densities at which they are optimal and the
proposition (lemma) where their asymptotic properties are described}%
\label{summary}
\begin{tabular*}{\tablewidth}{@{\extracolsep{\fill}}llllll@{}}
\hline
Test statistic& Reference  &Location $\theta$ &Validity (comments)&Optimality (comments)&Asymptotic  properties\\
\hline
$S\n_1(\theta)$ or $m_3^{(n)}(\theta)$ &(\ref{classicaltheta})
&Specified &Any $g_1$ with finite &Suboptimal at the Gaussian & Proposition
\ref{testGbart2} \\
(classical)&&&6th-order moments && \\
$S\n_2$ or $b_1^{(n)} $ & (\ref{classical}) &Unspecified &Any $g_1$
with finite & At the Gaussian, uniformly in $\theta, \sigma$ & Proposition
\ref{prop3.4} \\
(classical)&&&6th-order moments &($={T}^{(n)\dagger}({\bar X}\n)
+\mathrm{o}_{\mathrm{P}}(1)$) &(\ref{S=T}) \\
$T\n_{f_1}(\theta, {\hat\sigma}_{\#} )$ & (\ref{stat1art1})
&Specified & Specified $f_1\in\mathcal{F}_1$ & At $f_\sigma$,
uniformly in $\sigma$ & Proposition \ref{test1art1} \\[3pt]
$T\n_{f_1}({\hat\theta}_{\#} , {\hat\sigma}_{\#} )$ & (\ref{stat1art1})&Unspecified &Specified $f_1\in\mathcal{F}_1$ & At
$f_\sigma$, uniformly in $\theta, \sigma$ & Proposition \ref{test1art1} \\[3pt]
$\hat{T}\n_{f_1}({\theta} , {\hat\sigma}_{\#} )$ & (\ref
{stat2aart1})&Specified & Any $g_1 \in\mathcal{F}^*_{f_1}$ (\ref
{casuffit}) & At $f_\sigma$, uniformly in $\sigma$ (bad behavior &
Lemma \ref{test2art1}  \\
&&&&under local alternatives and $g_1\neq f_1$)&\\
${T}^{(n)\circ}_{f_1}({\theta} , {\hat\sigma}_{\#} )$ &(\ref
{stat2bart1}) &Specified &Any $g_1\in\mathcal{F}_1^{\circ}$
(\ref{casuffit2}) & At $f_\sigma$, uniformly in $\sigma$ (same as
$\hat{T}\n_{f_1}$,& Proposition \ref{test2bart1} \\
&&&& except for its bad behavior under & \\
&&&& local alternatives and $g_1\neq f_1$)&\\
${T}^{(n)\circ}_{f_1}({\hat\theta}_{\#} , {\hat\sigma}_{\#}
)$ &(\ref{stat2bart1}) &Unspecified &Any $g_1\in\mathcal{F}_1^{\circ}$ (\ref{casuffit2}) & At $f_\sigma$, uniformly in $\theta,
\sigma$ & Proposition \ref{test3art1} \\[3pt]
${T}^{(n)\dagger}(\theta) $&(\ref{statGbart1}) &Specified &Any $g_1
\in\mathcal{F}_{\phi_1}$ (implies & At the Gaussian, uniformly in $
\sigma$ & Proposition \ref{prop3.4}\\
(pseudo-Gaussian)&&&finite 6-order moments)&&
\\
${T}^{(n)\dagger}({\hat\theta} ) $&(\ref{statGbart1}) &Unspecified
&Any $g_1 \in\mathcal{F}_{\phi_1}$ (implies & At the Gaussian,
uniformly in $\theta, \sigma$ & Proposition \ref{testGbart2} \\
(pseudo-Gaussian)&&&finite 6-order moments)&&
\\
${T}^{(n)\circ }_\mathcal{L}(\theta) $&(\ref{TLaplace})
&Specified &Any $g_1 \in\mathcal{F}_{f_\mathcal{L}}$ (implies & At
the double exponential, & Proposition \ref{prop3.6} \\
(Laplace)&&&finite 4-order moments)& uniformly in $ \sigma$&
\\
${T}^{(n)\circ }_\mathcal{L}({\hat\theta} ) $&(\ref
{TLaplace}) &Unspecified &Any $g_1 \in\mathcal{F}_{f_\mathcal{L}}$
(implies & At the double exponential, & Proposition \ref{prop3.6} \\
(Laplace)&&&finite 4-order moments)&uniformly in $\theta, \sigma$ &
\\
\hline
\end{tabular*}
\end{sidewaystable}

These tests naturally extend into nonparametric rank-based ones.
The hypothesis of symmetry indeed enjoys strong group invariance
features. The null hypothesis $\mathcal{H}\n_\theta$ of symmetry
with respect to $\theta$ is generated by the group $\mathcal{G}\n
_\theta ,\sirc$ of all transformations~$\scriptscriptstyle\mathcal
{G}_h$ of $\R^n$ such that ${\scriptscriptstyle\mathcal
{G}}_h(x_1,\ldots, x_n):=(h(x_1),\ldots, h(x_n))$, where $\lim
_{x\rightarrow\pm\infty}h(x)=\pm\infty$, and $x\mapsto h(x)$ is
continuous, monotone increasing and skew-symmetric with respect to
$\theta$ (that is, satisfies $h(\theta-z)-\theta=-(h(\theta+
z)-\theta)$). A maximal invariant for that group is known to be the
vector $(s_1(\theta),\ldots, s_n(\theta))$, along with the vector $
(R\n_{+,1}(\theta),\ldots, R\n_{+,n}(\theta))$, where~$s_i(\theta
)$ is the sign of $X_i-\theta$ and $R\n_{+,i}(\theta)$ the rank of
$\vert X_i-\theta\vert$ among $\vert X_1-\theta\vert, \ldots,
\vert X_n-\theta\vert$. General results on semi-parametric efficiency
\cite{HW2003} indicate that, in such a context, the
expectation of the central sequence $\bolds{\Delta}\n_{f_1}(\varthetab
)$ conditional on those \textit{signed ranks} yields a~version of the
semi-parametrically efficient (at $f_1$ and $\varthetab$) central sequence.

That rank-based approach is adopted in a companion paper \cite{CHP2010}. For instance, the rank-based counterpart of the
specified-$\theta$ test statistic of Proposition \ref{testGbart1}(ii)
is the (strictly distribution-free, irrespective of any moment
assumptions) van der Waerden test based~on
\[
\utT\n_{\mathrm{vdW}}(\theta) :=\frac{1}{\sqrt {n \utg \n(\phi_1)}} \sum
_{i=1}^n s_i (\theta)\Phi^{-1} \!\biggl(\frac{n+1+R_{+, i}^{(n)}(\theta
)}{2(n+1)} \biggr)\! \biggl( \!\biggl(\Phi^{-1} \!\biggl(\frac{n+1+R_{+, i}^{(n)}(\theta)}{2(n+1)}
\biggr)\! \biggr)^2 \!-\! 3 \biggr),
\]
where $\utg \n(\phi_1):=n^{-1}\sum_{r=1}^n\Phi^{-1} (\frac
{n+1+r}{2(n+1)} ) ( (\Phi^{-1} (\frac{n+1+r}{2(n+1)} ) )^2 - 3 )^2$
and $\Phi$ stands for the standard normal distribution function. The
unspecified-$\theta$ case under such approach, however, is
considerably more delicate.

\begin{appendix}
\section*{Appendix} \label{appendixart1}
\subsection{\texorpdfstring{Proof of Proposition \protect\ref{lanart1}}{Proof of Proposition 2.1}}
The proof relies on \cite{Swensen1985}, Lemma 1, which involves a set of
six jointly sufficient conditions. Most of them readily follow from the
form of local likelihoods and are left to the reader. The most delicate
one is the quadratic mean differentiability of $(\theta, \sigma, \xi
) \mapsto g_{\theta, \sigma, \xi;f_1}^{1/2}(x)$,\vspace*{1pt} which we establish
in the following lemma, where $g_{\theta, \sigma, \xi;f_1} (x)$ is
the density defined in~(\ref{modelart1}).
\begin{Lem}\label{lemma1art1}
Let $f_1\in\mathcal{F}_1$, $\theta\in\R$, $\sigma\in\R^+_0$ and
$\xi\in\R$. Define
\begin{eqnarray*}
g _{\theta, \sigma, \xi; f_1}(x)&:=&\sigma^{-1} f_1 \biggl(\frac
{x-\theta}{\sigma} \biggr) -\frac{ \xi}{\sigma} {\dot{f}}_1 \biggl(\frac
{x-\theta}{\sigma} \biggr) \biggl( \biggl(\frac{x-\theta}{\sigma} \biggr)^2-\kappa(f_1)
\biggr)I[\vert x -\theta\vert\leq\sigma\vert z^* \vert]
\\
&&{} - \operatorname{sign}(\xi) \sigma^{-1} f_1 \biggl(\frac{x-\theta}{\sigma} \biggr) \{
I[ x-\theta>\operatorname{sign}(-\xi)\sigma|z^*| ]\\
&&\hphantom{{} - \operatorname{sign}(\xi) \sigma^{-1} f_1 \biggl(\frac{x-\theta}{\sigma} \biggr) \{}
{}- I[ x-\theta<\operatorname{sign}(\xi)\sigma|z^*| ] \},\nonumber\\
D_{\theta}g^{1/2}_{\theta, \sigma, 0; f_1}(x)&:=&\frac{1}{2} \sigma
^{-3/2} f_1^{1/2} \biggl(\frac{x-\theta}{\sigma} \biggr) \phi_{f_1} \biggl(\frac
{x-\theta}{\sigma} \biggr), \\
D_{\sigma}g^{1/2}_{\theta, \sigma, 0; f_1}(x)&:=&\frac{1}{2} \sigma
^{-3/2} f_1^{1/2} \biggl(\frac{x-\theta}{\sigma} \biggr) \biggl( \biggl(\frac{x-\theta
}{\sigma} \biggr) \phi_{f_1} \biggl(\frac{x-\theta}{\sigma} \biggr) -1 \biggr) 
\end{eqnarray*}
and
\begin{eqnarray*}
D_{\xi}g^{1/2}_{\theta, \sigma, \xi; f_1}(x)\vert_{\xi
=0}&:=&\frac{1}{2} \sigma^{-1/2} f_1^{1/2} \biggl(\frac{x-\theta}{\sigma
} \biggr) \phi_{f_1} \biggl(\frac{x-\theta}{\sigma} \biggr) \biggl( \biggl(\frac{x-\theta
}{\sigma} \biggr)^2-\kappa(f_1) \biggr).
\end{eqnarray*}
Then,
as $r,s$ and $t\to0$,
\begin{eqnarray*}
\mbox{\textup{(i)}}&\quad\hspace*{-4pt}&\int\{g^{1/2}_{\theta+t, \sigma+s, r;
f_1}(x)-g^{1/2}_{\theta+t, \sigma+s, 0; f_1}(x)- r D_{\xi
}g^{1/2}_{\theta+t, \sigma+s, \xi; f_1}(x)\vert_{\xi=0}\}^2 \,\mathrm{d}x =\mathrm{o}
(r^2).\\
\mbox{\textup{(ii)}}&\quad\hspace*{-4pt}&\int \left\{g^{1/2}_{\theta+t, \sigma+s, 0;
f_1}(x)-g^{1/2}_{\theta, \sigma, 0; f_1}(x)- \pmatrix{ t \cr s
}' \pmatrix{ D_{\theta}g^{1/2}_{\theta, \sigma, 0; f_1}(x) \cr
D_{\sigma}g^{1/2}_{\theta, \sigma, 0; f_1}(x)} \right\}^2 \,\mathrm{d}x =\mathrm{o} \left( \left\Vert \pmatrix{ t \cr s} \right\Vert^2
\right).\\
\mbox{\textup{(iii)}} &\quad\hspace*{-4pt}&\int \{ ( D_{\xi}g^{1/2}_{\theta+t, \sigma+s,
\xi; f_1}(x)\vert_{\xi=0} - D_{\xi}g^{1/2}_{\theta, \sigma, \xi;
f_1}(x)\vert_{\xi=0} ) \}^2 \,\mathrm{d}x =\mathrm{o} (1).\\
\mbox{\textup{(iv)}}&\quad\hspace*{-4pt}&\int \left\{g^{1/2}_{\theta+t, \sigma+s, r;
f_1}(x)-g^{1/2}_{\theta, \sigma, 0; f_1}(x)- \pmatrix{ t \cr s \cr r
}' \pmatrix{ D_{\theta}g^{1/2}_{\theta, \sigma, 0; f_1}(x) \cr
D_{\sigma}g^{1/2}_{\theta, \sigma, 0; f_1}(x) \cr D_{\xi
}g^{1/2}_{\theta, \sigma, \xi; f_1}(x)\vert_{\xi=0}
} \right\}^2 \,\mathrm{d}x\\
&&\quad =\mathrm{o} \left( \left\Vert \pmatrix{ t \cr s \cr r
} \right\Vert^2 \right).
\end{eqnarray*}
\end{Lem}

\begin{pf}
(i) Decompose
$\int \{g^{1/2}_{\theta+t, \sigma+s, r;
f_1}(x)-g^{1/2}_{\theta+t, \sigma+s, 0; f_1}(x)- r D_{\xi
}g^{1/2}_{\theta+t, \sigma+s, \xi; f_1}(x) \}^2 \,\mathrm{d}x$ into
$ a_1 +2a_2$,
where
\begin{eqnarray*}
a_1&=&\int_{|u|<|z^*|} \biggl\{ \biggl(\frac{1}{\sigma+s} f_1(u) \biggr)^{1/2} \bigl[ 1+r
\phi_{f_1} (u) \bigl(u^2-\kappa(f_1)\bigr) \bigr]^{1/2}- \biggl(\frac{1}{\sigma+s}
f_1(u) \biggr)^{1/2} \\
&&\hphantom{\int_{|u|<|z^*|} \biggl\{}
{} - \frac{r}{2} (\sigma+s)^{-1/2} \frac{\dot
{f}_1(u)}{f_1^{1/2}(u)} \bigl(u^2-\kappa(f_1)\bigr) \biggr\}^2 (\sigma+s)\,\mathrm{d}u
\end{eqnarray*}
and
\begin{eqnarray*}
a_2&=&\int_{u>|z^*|} \biggl\{ \bigl((\sigma+s)^{-1} f_1(u) \bigr)^{1/2} - \frac
{r}{2} (\sigma+s)^{-1/2} \frac{\dot{f}_1(u)}{f_1^{1/2}(u)}
\bigl(u^2-\kappa(f_1)\bigr) \biggr\}^2 (\sigma+s)\,\mathrm{d}u.
\end{eqnarray*}
Since, for $|x|<1$, $(1+x)^{1/2}=1+\frac{x}{2}(1+\lambda x)^{-1/2}$
for some $\lambda\in(0, 1)$, one easily obtains that
\begin{eqnarray*}
a_1&=&\frac{r^2}{4}\int_{|u|<|z^*|} \biggl\{ (\sigma+s)^{-1/2} \frac{\dot
{f}_{1}(u)}{f_1^{1/2}(u)} \bigl(u^2 -\kappa(f_1)\bigr)\\
&&\hphantom{\frac{r^2}{4}\int_{|u|<|z^*|} \biggl\{}
{}\times \bigl(\bigl(1+ \lambda r \phi
_{f_1} (u) \bigl(u^2 -\kappa(f_1)\bigr)\bigr)^{-1/2} -1\bigr) \biggr\}^2
(\sigma+s)\,\mathrm{d}u .
\end{eqnarray*}
For $|u|<1$, one has $(1-(1+\lambda u)^{-1/2})^2 \leq2 \frac
{2-\lambda}{1-\lambda}$, and the integrand is bounded by
\[
2 \frac{2-\lambda}{1-\lambda} \bigl(u^2-\kappa(f_1)\bigr)^2 \biggl(\frac{\dot
{f}_{1}^{1/2}(u)}{f_1^{1/2}(u)} \biggr)^2 ,
\]
which is square-integrable; Lebesgue's dominated convergence theorem
thus implies that~$a_1$ is $\mathrm{o}(r^2)$. Turning to $a_2$, we have that
$a_2 \leq C((\sigma+s)^{-1}a_{21}+a_{22})$, where
\[
a_{21}:=\int_{u>|z^*|} f_1(u) \,\mathrm{d}u \quad \mbox{and}\quad
a_{22}:=\frac{r^2}{4} \int_{u>|z^*|} \biggl(\frac{\dot
{f}_1(u)}{f_1^{1/2}(u)} \biggr)^2 \bigl(u^2-\kappa(f_1)\bigr)^2 \,\mathrm{d}u.
\]
The definition of $\mathcal{F}_1$ implies that $a_{21}=\mathrm{O}((z^*)^{-\beta
})$, hence that $a_{21}=\mathrm{o}(r^2)$ if $r (z^*)^{\beta/2} \rightarrow
\infty$ as $r \rightarrow0$. This latter condition holds, since $\phi
_{f_1}(z)=\mathrm{o}(z^{\beta/2-2})$ and since the definition of~$z^*$ entails that
\[
-1=r(z^*)^{\beta/2}\frac{\phi_{f_1}(z^*)}{z^{*(\beta/2-2)}}\frac
{(z^{*2}-\kappa(f_1))}{z^{*2}}.
\]
An application of Lebesgue's dominated convergence theorem again
yields $a_{22}=\mathrm{o}(r^2)$.

(ii) This is a particular case of Lemma A.1 in \cite{HP2006} (here in a simpler univariate context).

(iii) The fact that $D_{\xi}g^{1/2}_{\theta, \sigma, \xi;
f_1}(x)\vert_{\xi=0}$ is square-integrable implies that
\[
\|D_{\xi}g^{1/2}_{\theta+t, \sigma, \xi; f_1}(x)\vert_{\xi
=0}-D_{\xi}g^{1/2}_{\theta, \sigma, \xi; f_1}(x)\vert_{\xi=0}\|_{L^2}=\mathrm{o}(1)
\]
as $t$ tends to zero. Define $f_{1;\exp}(x):=f_1(\mathrm{e}^x)$
and $(f_{1;\exp}^{1/2}(x))\pr:=\frac
{1}{2}f_{1}^{-1/2}(\mathrm{e}^x)\dot{f}_{1}(\mathrm{e}^x)\mathrm{e}^x $. For the perturbation of
$\sigma$, we have
\begin{eqnarray*}
&&\int \bigl\{ \bigl( D_{\xi}g^{1/2}_{\theta, \sigma+s, \xi; f_1}(x)\vert
_{\xi=0} - D_{\xi}g^{1/2}_{\theta, \sigma, \xi; f_1}(x)\vert_{\xi
=0} \bigr) \bigr\}^2 \,\mathrm{d}x\\
&&\quad  =  2\sigma\int_0^{\infty} \biggl|\sigma^{-1/2} z \biggl(\biggl(1+\frac{s}{\sigma
}\biggr)^{-3/2} (f_{1;\exp}^{1/2})\pr\biggl(\ln(z)-\ln\biggl(1+\frac
{s}{\sigma}\biggr)\biggr)
- (f_{1;\exp}^{1/2})\pr(\ln(z)) \biggr)\\
&&\quad \hphantom{2\sigma\int_0^{\infty} \biggl|{}-}
-\sigma^{-1/2} z^{-1} \kappa(f_1)\\
&&\qquad \hphantom{2\sigma\int_0^{\infty} \biggl|}
\hspace{7pt}{}\times \biggl(\biggl(1+\frac{s}{\sigma}\biggr)^{1/2}
(f_{1;\exp}^{1/2})\pr\biggl(\ln(z)-\ln\biggl(1+\frac{s}{\sigma}\biggr)\biggr)
+(f_{1;\exp}^{1/2})\pr(\ln(z)) \biggr) \biggr|^2 \,\mathrm{d}z\\
& &\quad  \leq C(c_1+c_2),
\end{eqnarray*}
where
\[
c_1=\int_{-\infty}^{\infty} \biggl(\mathrm{e}^{(3/{2})(u-\ln(1+
{s}/{\sigma}))} (f_{1;\exp}^{1/2})\pr\biggl(u-\ln\biggl(1+\frac
{s}{\sigma}\biggr)\biggr)
- \mathrm{e}^{({3}/{2})u}(f_{1;\exp}^{1/2})\pr(u) \biggr)^2 \,\mathrm{d}u
\]
and
\[
c_2=\int_{-\infty}^{\infty} \biggl(\mathrm{e}^{-(1/{2})(u-\ln(1+
{s}/{\sigma}))} (f_{1;\exp}^{1/2})\pr\biggl(u-\ln\biggl(1+\frac
{s}{\sigma}\biggr)\biggr)
- \mathrm{e}^{-(1/{2})u}(f_{1;\exp}^{1/2})\pr(u) \biggr)^2 \,\mathrm{d}u.
\]
Now, both $\mathrm{e}^{-(1/{2})u}(f_{1;\exp}^{1/2})\pr(u)$
and $\mathrm{e}^{({3}/{2})u}(f_{1;\exp}^{1/2})\pr(u)$ are
square-integrable since $f_1 \in\mathcal{F}_1$. Therefore, quadratic
mean continuity implies that $c_1$ and $c_2$ are $\mathrm{o}(1)$ as $s
\rightarrow0$.

(iv) The left-hand side in (iv) is bounded by
$C(b_1+b_2+b_3)$, where
\begin{eqnarray*}
b_1 &=&\int \{g^{1/2}_{\theta+t, \sigma+s, r;
f_1}(x)-g^{1/2}_{\theta+t, \sigma+s, 0; f_1}(x)
- r D_{\xi}g^{1/2}_{\theta+t, \sigma+s, \xi; f_1}(x) \}^2 \,\mathrm{d}x,
\\
b_2&=&\int \left\{g^{1/2}_{\theta+t, \sigma+s, 0; f_1}(x)-g^{1/2}_{\theta,
\sigma, 0; f_1}(x)- ( t ,\ s ) \pmatrix{
D_{\theta}g^{1/2}_{\theta, \sigma, 0; f_1}(x) \cr
D_{\sigma}g^{1/2}_{\theta, \sigma, 0; f_1}(x)
} \right\}^2 \,\mathrm{d}x
\end{eqnarray*}
and
\[
b_3=\int \bigl\{r \bigl( D_{\xi}g^{1/2}_{\theta+t, \sigma+s, \xi; f_1}(x) -
D_{\xi}g^{1/2}_{\theta, \sigma, \xi; f_1}(x) \bigr) \bigr\}^2 \,\mathrm{d}x.
\]
The result then follows from (i)--(iii).
\end{pf}

\subsection{Asymptotic linearity}
\subsubsection[Asymptotic linearity of Delta (n) f1;3]{Asymptotic linearity of $\Delta\n_{f_1;3}$}
The asymptotic linearity of $\Delta\n_{f_1;3}$ is required in the
construction of the optimal parametric test of Section \ref{322}. Note
that the proof below needs uniform local asymptotic normality in
$\theta$ and $\sigma$ only.

\begin{Prop} \label{lin1art1}
Let $f_1 \in\mathcal{F}_1$ and $g_1 \in\mathcal{F}_{f_1}$. Then,
under ${\mathrm{P}}^{(n)}_{\theta, \sigma, 0; g_1}$, as $n \rightarrow
\infty$,
\begin{enumerate}[(ii)]
\item[(i)] $
\Delta_{f_1;3}^{(n)}(\theta+n^{-1/2} t, \sigma, 0) = \bolds{\Delta
}_{f_1;3}^{(n)}(\theta,\sigma, 0)- t \sigma^{-1} (\mathcal
{J}_{g_1}(f_1, g_1) - \kappa(f_1)\mathcal{I}_{g_1}(f_1, g_1) )
+\mathrm{o}_{\mathrm{P}}(1)$ for all $t \in\mathbb R$.
\item[(ii)] 
$\bolds{\Delta}_{f_1;3}^{(n)}(\theta, \sigma+n^{-1/2} s, 0) = \bolds{\Delta}_{f_1;3}^{(n)}(\theta,\sigma, 0)
+\mathrm{o}_{\mathrm{P}}(1)$ for all $s \in\mathbb R$.
\end{enumerate}
\end{Prop}

\begin{pf}
Define
\[
D_{f_1;3}^{(n)}(\theta, \sigma):=n^{-1/2} \sum_{i=1}^n \phi
_{f_1}\bigl(Z_i^{(n)}(\theta, \sigma)\bigr) \bigl(Z_i^{(n)}(\theta, \sigma)\bigr)^2.
\]
Letting $K_{f_1}(u):=\phi_{f_1}(G_{1+}^{-1}(u))$ $(G_{1+}^{-1}(u))^2$ with $G_{1+}(z):= (2G_1(z)-1)I[z\geq 0]$,
note that
\[
D_{f_1;3}^{(n)}(\theta, \sigma, 0)=n^{-1/2} \sum_{i=1}^n s_i (\theta
)K_{f_1}\bigl(G_{1+}\bigl(\bigl|Z_i^{(n)}(\theta, \sigma)\bigr|\bigr)\bigr).
\]
Writing $\theta\n$ for $\theta+n^{-1/2}t$, $Z_{i}^{(n)}$ for
$Z_i^{(n)}(\theta, \sigma)$, $s_i^{(n)}$ for $\mbox
{sign}(Z_{i}^{(n)})$, $Z_{i;n}^{(n)}$ for $Z_i^{(n)}(\theta\n, \sigma
)$ and $s_{i;n}^{(n)}$ for $\mbox{sign}(Z_{i;n}^{(n)})$, let us show that
\[
n^{-1/2} \sum_{i=1}^n s_{i;n}^{(n)}
K_{f_1}\bigl(G_{1+}\bigl(\bigl|Z_{i;n}^{(n)}\bigr|\bigr)\bigr)-n^{-1/2} \sum_{i=1}^n s_i^{(n)}
K_{f_1}\bigl(G_{1+}\bigl(\bigl|Z_{i}^{(n)}\bigr|\bigr)\bigr)+t \sigma^{-1} \mathcal
{J}_{g_1}(f_1,g_1)
\]
is $\mathrm{o}_{\mathrm{P}}(1)$; the proofs of (ii) and that of
\[
n^{-1/2} \sum_{i=1}^n \phi_{f_1}\bigl(Z_{i;n}^{(n)}\bigr) = n^{-1/2} \sum
_{i=1}^n \phi_{f_1}\bigl(Z_i^{(n)}\bigr)- t \sigma^{-1} \mathcal{I}_{g_1}(f_1, g_1)
+\mathrm{o}_{\mathrm{P}}(1)
\]
follow along the same lines and are therefore left to the reader.
Let
\begin{eqnarray*}
K_{f_1}^{(l)}(u)&:=& K_{f_1}\biggl(\frac{2}{l}\biggr)l\biggl(u-\frac{1}{l}\biggr) I\biggl[\frac
{1}{l} < u \leq\frac{2}{l}\biggr] +K_{f_1}(u) I\biggl[\frac{2}{l} < u \leq
1-\frac{2}{l}\biggr]\\
&&{}+ K_{f_1}\biggl(1-\frac{2}{l}\biggr)l\biggl(\biggl(1-\frac{1}{l}\biggr)-u\biggr) I\biggl[1-\frac{2}{l} < u
\leq1-\frac{1}{l}\biggr].
\end{eqnarray*}
Continuity of $u\mapsto K_{f_1}(u)$ implies continuity of $u\mapsto
K_{f_1}^{(l)}(u)$ on the interval $]0,1[$. Moreover, since this
function is compactly supported, it is bounded, for any (sufficiently
large) $l \in\mathbb N_0$, by the monotone increasing function
$u\mapsto K_{f_1}(u)$.
Let $\mathrm{E}_0$ denote expectation under ${\mathrm{P}}^{(n)}_{\theta
,\sigma,0; g_1}$. One easily shows that $D_{f_1;3}^{(n)}(\theta,
\sigma)$ decomposes into
$D_1^{(n,m)}+D_2^{(n,m)}-R_1^{(n,m)}+R_2^{(n,m)}+ R_3^{(m)}$ where,
defining
$\mathcal{J}_{g_1}^{(l)}(f_1, g_1):=\int_0^1 K_{f_1}^{(l)} (u) \phi
_{g_1} ( G_1^{-1}(u)) \,\mathrm{d}u$,
\begin{eqnarray*}
D_1^{(n,l)}&=&n^{-1/2} \sum_{i=1}^n s_{i;n} K_{f_1}^{(l)}
\bigl(G_{1+}\bigl(\bigl|Z_{i;n}^{(n)}\bigr|\bigr)\bigr)- n^{-1/2} \sum_{i=1}^n s_{i} K_{f_1}^{(l)}
\bigl(G_{1+}\bigl(\bigl|Z_{i}^{(n)}\bigr|\bigr)\bigr)\\
&&{} -n^{1/2} \mathrm{E}_0\bigl[s_{i;n}K_{f_1}^{(l)}
\bigl(G_{1+}\bigl(\bigl|Z_{i;n}^{(n)}\bigr|\bigr)\bigr)\bigr], \\
D_2^{(n,l)}&=&n^{1/2} \mathrm{E}_0\bigl[s_{i;n}K_{f_1}^{(l)}
\bigl(G_{1+}\bigl(\bigl|Z_{i;n}^{(n)}\bigr|\bigr)\bigr)\bigr] + t \sigma^{-1} \mathcal{J}_{g_1}^{(l)}(f_1,
g_1),\\
R_1^{(n,l)}&=& n^{-1/2} \sum_{i=1}^{n} s_{i} \bigl(K_{f_1}
\bigl(G_{1+}\bigl(\bigl|Z_{i}^{(n)}\bigr|\bigr)\bigr) - K_{f_1}^{(l)} \bigl(G_{1+}\bigl(\bigl|Z_{i}^{(n)}\bigr|\bigr)\bigr)\bigr),\\
R_2^{(n,l)}&=& n^{-1/2} \sum_{i=1}^{n}s_{i;n} \bigl(K_{f_1}
\bigl(G_{1+}\bigl(\bigl|Z_{i;n}^{(n)}\bigr|\bigr)\bigr) - K_{f_1}^{(l)} \bigl(G_{1+}\bigl(\bigl|Z_{i;n}^{(n)}\bigr|\bigr)\bigr)\bigr)
\end{eqnarray*}
and
%
\[
R_3^{(l)}= t \sigma^{-1} \bigl(\mathcal{J}_{g_1}(f_1, g_1)-\mathcal
{J}_{g_1}^{(l)}(f_1, g_1)\bigr).
\]
In order to conclude, we prove that $D_1^{(n,l)}$ and $D_2^{(n,l)}$ are
$\mathrm{o}_{\mathrm{P}}(1)$ under ${\mathrm{P}}^{(n)}_{\theta,\sigma,0; g_1}$,
as \mbox{$n\rightarrow\infty$}, for fixed $l$, and that $R_1^{(n,l)}$,
$R_2^{(n,l)}$ and $R_3^{(l)}$ are $\mathrm{o}_{\mathrm{P}}(1)$ under the same
sequence of hypotheses, as $l \rightarrow\infty$, uniformly in $n$.
For the sake of convenience, those three results are treated separately
(Lemmas \ref{lemme1art1}--\ref{lemme3art1}).

\begin{Lem} \label{lemme1art1}
For any fixed $l$,
$D_1^{(n,l)}=\mathrm{o}_{\mathrm{P}}(1)$ as $n \rightarrow\infty$, under ${\mathrm{P}}^{(n)}_{\theta,\sigma,0; g_1}$.
\end{Lem}
\begin{Lem} \label{lemme2art1}
For any fixed $l$,
$D_2^{(n,l)}=\mathrm{o}_{\mathrm{P}}(1)$ as $n \rightarrow\infty$, under ${\mathrm{P}}^{(n)}_{\theta,\sigma,0; g_1}$.
\end{Lem}
\begin{Lem} \label{lemme3art1}
\textup{(i)} Under ${\mathrm{P}}^{(n)}_{\theta,\sigma,0; g_1}$, $R_1^{(n,l)}=\mathrm{o}_{\mathrm{P}}(1)$ as $l \rightarrow\infty$, uniformly in $n$.
{\smallskipamount=0pt
\begin{longlist}[(iii)]
\item[(ii)] $R_2^{(n,l)}=\mathrm{o}_{\mathrm{P}}(1)$ as $l \rightarrow\infty$,
under ${\mathrm{P}}^{(n)}_{\theta,\sigma,0; g_1}$ (for $n$ sufficiently
large), uniformly in $n$.
\item[(iii)] $R_3^{(l)}$ is $\mathrm{o}(1)$ as $l \rightarrow\infty$.
\end{longlist}
}
\end{Lem}

\begin{pf*}{Proof of Lemma \ref{lemme1art1}}
Consider the i.i.d. variables
\[
T_i^{(n,l)}:=s_{i;n}
K_{f_1}^{(l)}\bigl(G_{1+}\bigl(\bigl|Z_{i;n}^{(n)}\bigr|\bigr)\bigr)-s_{i}
K_{f_1}^{(l)}\bigl(G_{1+}\bigl(\bigl|Z_{i}^{(n)}\bigr|\bigr)\bigr).
\]
One easily verifies that
$D_1^{(n,l)}=n^{-1/2} \sum_{i=1}^n (T_i^{(n,l)}-\mathrm{E}_0[T_i^{(n,l)}])$. Writing $\operatorname{Var}_0$ for variances under ${\mathrm{P}}^{(n)}_{\theta,\sigma,0; g_1}$, we have that
\begin{eqnarray*}
\mathrm{E}_0\bigl(D_{1}^{(n,l)}\bigr) &\leq&n^{-1}\mathrm{E}_0\Biggl[\Biggl(\sum
_{i=1}^n\bigl(T_i^{(n,l)}-\mathrm{E}_0\bigl[T_i^{(n,l)}\bigr]\bigr)\Biggr)^2\Biggr]\\
&\leq& n^{-1}\operatorname{Var}_0\Biggl[\sum_{i=1}^n\bigl(T_i^{(n,l)}-\mathrm{E}_0\bigl[T_i^{(n,l)}\bigr]\bigr)\Biggr]
= \operatorname{Var}_0\bigl[T_i^{(n,l)}\bigr]
\leq\mathrm{E}_0\bigl[\bigl(T_i^{(n,l)}\bigr)^2\bigr],
\end{eqnarray*}
and it only remains to show that
\[
\mathrm{E}_0\bigl[\bigl(T_i^{(n,l)}\bigr)^2\bigr]=\mathrm{E}_0\bigl[\bigl(s_{i;n}
K_{f_1}^{(l)}\bigl(G_{1+}\bigl(\bigl|Z_{i;n}^{(n)}\bigr|\bigr)\bigr)-s_{i}
K_{f_1}^{(l)}\bigl(G_{1+}\bigl(\bigl|Z_{i}^{(n)}\bigr|\bigr)\bigr)\bigr)^2\bigr]=\mathrm{o}(1)
\]
as $n \rightarrow\infty$. Now,\vspace*{-3pt}
\begin{eqnarray*}
& & \bigl(s_{i;n} K_{f_1}^{(l)}\bigl(G_{1+}\bigl(\bigl|Z_{i;n}^{(n)}\bigr|\bigr)\bigr)-s_{i}
K_{f_1}^{(l)}\bigl(G_{1+}\bigl(\bigl|Z_{i}^{(n)}\bigr|\bigr)\bigr)\bigr)^2\\
&&\quad  =\bigl(s_{i;n}
K_{f_1}^{(l)}\bigl(G_{1+}\bigl(\bigl|Z_{i;n}^{(n)}\bigr|\bigr)\bigr)-s_{i;n}K_{f_1}^{(l)}\bigl(G_{1+}\bigl(\bigl|Z_{i}^{(n)}\bigr|\bigr)\bigr)
+s_{i;n}K_{f_1}^{(l)}\bigl(G_{1+}\bigl(\bigl|Z_{i}^{(n)}\bigr|\bigr)\bigr)\\
&&\qquad \hphantom{\bigl(}{}-s_{i}
K_{f_1}^{(l)}\bigl(G_{1+}\bigl(\bigl|Z_{i}^{(n)}\bigr|\bigr)\bigr)\bigr)^2\\
&&\quad  \leq2 \bigl(K_{f_1}^{(l)}\bigl(G_{1+}\bigl(\bigl|Z_{i;n}^{(n)}\bigr|\bigr)\bigr)-
K_{f_1}^{(l)}\bigl(G_{1+}\bigl(\bigl|Z_{i}^{(n)}\bigr|\bigr)\bigr)\bigr)^2+2\bigl(K_{f_1}^{(l)}\bigl(G_{1+}\bigl(\bigl|Z_{i}^{(n)}\bigr|\bigr)\bigr)\bigr)^2(s_{i;n}-s_{i})^2.
\vspace*{-3pt}
\end{eqnarray*}
Because $u\mapsto K_{f_1}^{(l)}(u)$ is continuous and
$|Z_{i;n}^{(n)}-Z_{i}^{(n)}|$ is $\mathrm{o}_{\mathrm{P}}(1)$,
$K_{f_1}^{(l)}(G_{1+}(|Z_{i;n}^{(n)}|))-
K_{f_1}^{(l)}(G_{1+}(|Z_{i}^{(n)}|))$ also is $\mathrm{o}_{\mathrm{P}}(1)$. Moreover,
since $K_{f_1}^{(l)}$ is bounded, this convergence to zero also holds
in quadratic mean. Similarly,
$K_{f_1}^{(l)}(G_{1+}(|Z_{i}^{(n)}|))(s_{i;n}-s_{i})=\mathrm{o}_{\mathrm{P}}(1)$ since
$K_{f_1}^{(l)}(G_{1+}(|Z_{i}^{(n)}|))$ is bounded and $|s_{i;n}-s_{i}|$
is $\mathrm{o}_{\mathrm{P}}(1)$. Finally, both $s_{i;n}$ and $s_{i}$ are bounded,
implying that this convergence to zero also holds in quadratic
mean.
\end{pf*}


\begin{pf*}{Proof of Lemma \ref{lemme2art1}}
Let $\DS{B_1^{(n,l)}:=n^{-1/2} \sum_{i=1}^n s_i
K_{f_1}^{(l)}(G_{1+}(|Z_{i}^{(n)}|))}$. As $n \rightarrow\infty$,
under ${\mathrm{P}}^{(n)}_{\theta,\sigma,0; g_1}$,
$\DS{B_1^{(n,l)}}$ is asymptotically ${ \mathcal{N}(0, \mbox
{E}[(K_{f_1}^{(l)}(U))^2])}$,
where $U$ stands for a random variable uniformly distributed over the
unit interval. 
Also, letting $\DS{B_2^{(n,l)}:=n^{-1/2} \sum_{i=1}^n s_{i;n}\times
K_{f_1}^{(l)}(G_{1+}|Z_{i;n}^{(n)}|)}$, it follows from ULAN that
$B_2^{(n,l)} - t \sigma^{-1} \mathcal{J}_{g_1}^{(l)}(f_1, g_1)$ is
asymptotically 
$\mathcal{N}(0, \mathrm{E}[(K_{f_1}^{(l)}(U))^2])$ as $n \rightarrow
\infty$, under ${\mathrm{P}}^{(n)}_{\theta,\sigma,0; g_1}$.
Since $D_1^{(n,l)}=B_2^{(n,l)}-B_1^{(n,l)}-\mbox
{E}_0[B_2^{(n,l)}]=\mathrm{o}_{\mathrm{P}}(1)$, we have that
$\DS{B_2^{(n,l)} - \mathrm{E}_0[B_2^{(n,l)}]}$ is asymptotically ${
\mathcal{N}(0, \mathrm{E}[(K_{f_1}^{(l)}(U))^2])}$
as $n \rightarrow\infty$, under ${\mathrm{P}}^{(n)}_{\theta,\sigma,0;
g_1}$. Therefore, still
as $n \rightarrow\infty$,
$D_2^{(n,l)}=\mathrm{E}_0[B_2^{(n,l)}]
- t \sigma^{-1} \mathcal{J}_{g_1}^{(l)}(f_1, g_1)=\mathrm{o}(1)$.
\end{pf*}

%


\begin{pf*}{Proof of Lemma \ref{lemme3art1}}
(i) We have that\vspace*{-3pt}
\begin{eqnarray*}
\mathrm{E}_0 \bigl[ \bigl(R_1^{(n,l)} \bigr)^2 \bigr]&\leq&C \mathrm{E}_0 \bigl[
\bigl(K_{f_1}\bigl(G_{1+}\bigl(\bigl|Z_{i}^{(n)}\bigr|\bigr)\bigr)-K_{f_1}^{(l)}\bigl(G_{1+}\bigl(\bigl|Z_{i}^{(n)}\bigr|\bigr)\bigr) \bigr)^2 \bigr]\\
&=& C \int_{-\infty}^{\infty} \bigl(K_{f_1}(u)-K_{f_1}^{(l)}(u) \bigr)^2 \,\mathrm{d}u;\vspace*{-3pt}
\end{eqnarray*}
for any $u \in\,]0,1[$, $K_{f_1}^{(l)}(u)$ converges to $K_{f_1}(u)$ and
the integrand is bounded (uniformly in~$l$) by $4 (K_{f_1}(u))^2$,
which is integrable on $]0,1[$. Lebesgue's dominated convergence
theorem thus implies that $\mathrm{E}_0[(R_1^{(n,l)})^2]=\mathrm{o}(1)$ as $l
\rightarrow\infty$, uniformly in $n$.

(ii) The claim here is the same as in (i), with
$Z_{i;n}^{(n)}$ replacing $Z_i^{(n)}$. Accordingly, (ii)~holds under
$\mathrm{P}^{(n)}_{\theta\n, \sigma, 0; g_1}$. That it also holds
under ${\mathrm{P}}^{(n)}_{\theta,\sigma,0; g_1}$ follows from Lemma 3.5
in \cite{J1969}.

(iii) Note that\vspace*{-3pt}
\begin{eqnarray*}
\bigl\vert\mathcal{J}_{g_1}(f_1, g_1)-\mathcal{J}_{g_1}^{(l)}(f_1, g_1)
\bigr\vert^2 &=& \biggl\vert\int_0^1 \phi_{g_1}(G_{1}^{-1}(u))
\bigl(K_{f_1}^{(l)}(u)-K_{f_1}(u)\bigr) \,\mathrm{d}u \biggr\vert^2
\\[-2pt]
&\leq& \mathcal{I}(g_1) \int_0^1 \bigl(\bigl(K_{f_1}^{(l)}(u)-K_{f_1}(u)\bigr)\bigr)^2 \,\mathrm{d}u,
\end{eqnarray*}
where the integrand is bounded by $4 (K_{f_1}(u))^2$, which is
square-integrable. Pointwise convergence of $K_{f_1}^{(l)}(u)$ to
$K_{f_1}(u)$ implies that $\mathcal{J}_{g_1}(f_1, g_1)-\mathcal
{J}_{g_1}^{(l)}(f_1, g_1)=\mathrm{o}(1)$ as $l \rightarrow\infty$. The result
follows.
\end{pf*}

\subsubsection[Substitution of Delta (n) f1;3 (kappa (n) (f1, theta, sigma), theta, sigma)]%
{Substitution of $\Delta_{f_1;3}^{(n)}(\kappa
^{(n){\circ} } (f_1, \hat{\theta}_\#,\hat{\sigma}_\# ),
\hat{\theta}_\#, \hat{\sigma}_\# )$ for $\Delta
_{f_1;3}^{(n)}(\kappa_{g_1}(f_1, g_1), \theta, {\sigma} )$.}

\begin{Lem}\label{subst}
Let $f_1\in\mathcal{F}_1^{\circ}$ and $g_1\in\mathcal
{F}_{f_1}^{\circ} $. Then, under ${\mathrm{P}}\n_{\theta, \sigma,
0; g_1}$:
\begin{longlist}[(iii)]
\item[(i)] $\Delta_{f_1;3}^{(n)}(\kappa_{g_1}(f_1, g_1), \hat
{\theta}_\#, \hat{\sigma}_\# )-\Delta_{f_1;3}^{(n)}(\kappa
_{g_1}(f_1, g_1), {\theta}, {\sigma} )=\mathrm{o}_{\mathrm{P}}(1)$.
\item[(ii)] $\Delta_{f_1;3}^{(n)}(\kappa^{(n){\circ} } (f_1,
\hat{\theta}_\#,\hat{\sigma}_\# ), \hat{\theta}_\#, \hat{\sigma
} _\#)-\Delta_{f_1;3}^{(n)}(\kappa_{g_1}(f_1, g_1), \hat{\theta}_\#
, \hat{\sigma}_\# )=\mathrm{o}_{\mathrm{P}}(1)$.
\item[(iii)] $\Delta_{f_1;3}^{(n)}(\kappa^{(n){\circ} }
(f_1, \hat{\theta}_\#,\hat{\sigma}_\# ), \hat{\theta}_\#, \hat
{\sigma}_\# )-\Delta_{f_1;3}^{(n)}(\kappa_{g_1}(f_1, g_1), \theta,
{\sigma} )=\mathrm{o}_{\mathrm{P}}(1)$.
\end{longlist}
\end{Lem}

\begin{pf}
Part (i) is a direct consequence of Proposition \ref
{lin1art1}.
The left-hand side in (ii) can be written as
%
\begin{equation}\label{299}
T_1\n\times T_2\n:= \bigl(\kappa^{(n){\circ} } (f_1, \hat{\theta
}_\#,\hat{\sigma}_\# )-\kappa_{g_1}(f_1, g_1)\bigr) \times n^{-1/2}\sum
_{i=1}^n \phi_{f_1}\bigl(Z_i^{(n)}(\hat{\theta}_\#,\hat{\sigma}_\#)\bigr)
\end{equation}
where $T_1\n$ is $\mathrm{o}_{\mathrm{P}}(1)$. Now, ULAN implies that
%
\begin{equation}\label{399}
T_2\n= n^{-1/2}\sum_{i=1}^n \phi_{f_1}\bigl(Z_i^{(n)}({\theta},{\sigma})\bigr)
+ ( \sigma^{-1} \mathcal{I}_{g_1}(f_1, g_1), 0 )n^{1/2} \left( \pmatrix{ \hat{\theta}_\# \cr \hat{\sigma }_\#
} - \pmatrix{ {\theta} \cr {\sigma }
} \right)+\mathrm{o}_{\mathrm{P}}(1),
\end{equation}
as $\ny$ under ${\mathrm{P}}\n_{\theta, \sigma, 0; g_1}$. Hence, the
central limit theorem and the root-$n$-consistency of $\hat{\theta}_\#
$ and $\hat{\sigma }_\#$ entail that (\ref{399}) is $\mathrm{O}_{\mathrm{P}}(1)$;
the result follows. As for (iii), it is a direct consequence of (i) and
(ii).
\end{pf}\noqed
\end{pf}
\end{appendix}

\section*{Acknowledgements}
Marc Hallin's research supported by the Sonderforschungsbereich
``Statistical modelling of nonlinear dynamic processes'' (SFB 823)
of the German Research Foundation (Deutsche Forschungsgemeinschaft) and
a Discovery Grant of the Australian Research Council. This work was
partly completed while visiting ORFE and the Bendheim Center at
Princeton University. Davy Paindaveine's research supported by a Mandat
d'Impulsion Scientifique of the Fonds National de la Recherche
Scientifique, Communaut\' e fran\c caise de Belgique.

\printhistory

\end{document}